\documentclass[10pt]{article}
\usepackage[utf8]{inputenc}

\usepackage{graphicx}%
\usepackage{multirow}%
\usepackage{amsmath,amssymb,amsfonts}%
\usepackage{amsthm}%
\usepackage{mathrsfs}%
\usepackage[title]{appendix}%
\usepackage{xcolor}%
\usepackage{textcomp}%
\usepackage{manyfoot}%
\usepackage{booktabs}%
\usepackage{algorithm}%
\usepackage{algorithmicx}%
\usepackage{algpseudocode}%
\usepackage{listings}%

\usepackage{amsthm,stmaryrd,mathtools}
\usepackage{graphicx}
\usepackage{multirow}
\usepackage{subcaption}

\renewcommand{\div}{\mathop{\rm div}\nolimits}

\graphicspath{ {./figs/} }

\theoremstyle{thmstyleone}%
\newtheorem{theorem}{Theorem}
\theoremstyle{thmstyletwo}%

\theoremstyle{thmstylethree}%

\begin{document}

\title{An adaptive two-grid preconditioner and linearly implicit scheme for shale gas transport in fractured porous media}

%
%
%
%
%

\author{
Maria Vasilyeva 
\thanks{Department of Mathematics and Statistics, Texas A\&M University  - Corpus Christi,   Corpus Christi, Texas, USA. Email: {\tt maria.vasilyeva@tamucc.edu}.} 
\and 
Ben S. Southworth  
\thanks{Theoretical Division, Los Alamos National Laboratory, NM, USA. Email: {\tt southworth@lanl.gov}.}
\and
Shubin Fu 
\thanks{Eastern Institute for Advanced Study, Eastern Institute of Technology, Zhejiang, China. Email: {\tt shubinfu@eias.ac.cn}.}
}

\maketitle

\begin{abstract}
We consider a nonlinear mixed-dimensional model for simulating gas transport in shale formation. The mathematical model consists of a coupled system of nonlinear equations, where flow within fractures is represented using a lower-dimensional representation.  For the numerical solution of the coupled transport problem, we construct an unstructured mesh that resolves lower dimensional fractures on the grid level and use the finite element approximation to build a discrete system. To construct an efficient scheme for the resulting nonlinear problem, we use an explicit-implicit method for time integration, where we carefully choose an additive partition of the nonlinear operators to separate the stiff linear component and integrate it implicitly to ensure the stability of the time integration. Next, we invert the linear partition of the operator by constructing an efficient two-grid preconditioner for shale gas transport in fractured porous media. We use a local pointwise smoother on the fine grid and carefully design an adaptive multiscale space for coarse grid approximation based on local generalized eigenvalue problems. We utilize an adaptive thresholding to automatically identify local dominant modes which correspond to the very small eigenvalues in local domains. We remark that such spatial features are automatically captured through our local spectral problems, and connect these to fracture information in the global formulation of the problem. Approximation properties of the local spectral space with convergence of the proposed two-grid algorithm are given. Numerical results are presented for two fracture distributions with 30 and 160 fractures, demonstrating iterative convergence independent of the contrast of fracture and porous matrix permeability.
\end{abstract}




\section{Introduction}

Numerical simulation of the complex processes in shale gas reservoirs is challenging due to the presence of multiple spatial and temporal scales. Moreover, shale formations are highly heterogeneous and contain a complex mixture of organic matter (kerogen), inorganic matter, and natural and hydraulic fractures \cite{akkutlu2012multiscale, wasaki2014permeability}. The transport mechanisms in shale gas reservoirs are complex and highly nonlinear, which include viscous flow and diffusion in addition to adsorption by the internal kerogen surfaces \cite{javadpour2007nanoscale, loucks2009morphology, freeman2011numerical, yao2013numerical, zhang2015triple, talonov2016numerical}. Moreover, the complex shale gas transport should be coupled to the transport in fractures. 
Fractured porous media is characterized by a complex flow direction associated with high-contrast properties \cite{akkutlu2018multiscale, vasilyeva2019multiscale, vasilyeva2019multiscale, lee2016pressure, hajibeygi2011hierarchical, wu2011multiple, yao2010discrete}. A common simulation is a multi-continuum formulation that describes problems with high-contrast coefficients, and separate flow based on heterogeneity \cite{barenblatt1960basic, arbogast1990derivation, vasilyeva2023efficient}. Dual porosity models describe flow in naturally fractured media and are common in reservoir simulation \cite{barenblatt1960basic, warren1963behavior}. Models with explicit fracture representation are usually used for large-scale hydraulic fractures, where a common approach for flow simulation is based on the mixed dimensional formulation \cite{martin2005modeling, d2012mixed, formaggia2014reduced, Quarteroni2008coupling, schwenck2015dimensionally}.  

Numerical simulations of the shale gas transport problem are challenging due to the nonlinear and multiscale nature of the underlying processes. In this work, we present the construction of an efficient numerical technique based on an additive partition of operators and an implicit-explicit time approximation scheme. Our work is motivated by the Explicit-Implicit-Null (EIN) scheme proposed in \cite{duchemin2014explicit}, where the authors proposed adding and subtracting linear terms to ensure the stability of the scheme. This additional term acts as a damping component and is treated implicitly, while other parts are treated explicitly. A similar approach was introduced earlier in  \cite{hou1994removing} for removing the stiffness from interfacial flows with surface tension by reformulating high-order, nonlinear, and non-local terms into linear. The Variable-Coefficient EIN method has been proposed for high-order diffusion and dispersion equations in \cite{tan2024high}. These methods add and subtract a linear highest derivative term, enabling stability through implicit-explicit discretization.
In \cite{wang2024partially}, multiscale, nonlinear, high-contrast diffusion problems are considered, and they utilize a partially explicit splitting scheme for multiscale approximation. 
The EIN scheme is closely related to linearly implicit methods \cite{hairer2010solving, calvo2001linearly}. 
Moreover, the EIN scheme can be interpreted as an additive scheme \cite{vabishchevich2013additive}, where the problem operator is represented as a sum of the linear and nonlinear parts with an implicit linear approximation of the stiff dynamics and an explicit approximation for the remaining nonlinear residual (implicit-explicit or ImEx scheme). ImEx scheme are used to integrate stiff processes in time without having to treat the entire operator implicitly, instead separating out non-stiff components for explicit integration, e.g. \cite{ascher1995implicit,ascher1997implicit,kennedy2019higher,kennedy2003additive,pareschi2005implicit}. A general analysis of additive multistep schemes is presented in \cite{vabishchevich2013additive} with a diverse range of applications of resulting splitting techniques. In our previous work, we used an ImEx scheme to construct a decoupled scheme for flow in fractured porous media \cite{vasilyeva2023efficient, vasilyeva2024implicit}. We proposed decoupling the resulting block matrix by explicitly approximating the coupling term. In \cite{vabishchevich2012explicit, afanas2013unconditionally}, we presented implicit-explicit schemes for convection-diffusion problems. 

The development of efficient solvers and preconditioner for nonlinear flow and transport in heterogeneous media has also been an active area of research, where the condition number can depend strongly on both the mesh size $h$ and the coefficient contrast. Traditional two-grid and multigrid preconditioning algorithms use coarse solvers based on linear or polynomial interpolations and may deteriorate in the presence of rapid small-scale oscillations or high aspect ratios. However, the design of robust preconditioners with contrast-independent performance tends to require special coarse spaces, e.g. \cite{g1, efendiev2010spectral,al2023efficient}. 
In \cite{aarnes2002multiscale}, the authors develop and analyze a class of nonoverlapping domain decomposition methods using multiscale coarse grid solvers using the Multiscale Finite Element Method (MsFEM). 
An efficient and robust multiscale preconditioner for large-scale incompressible flow in highly heterogeneous porous media is presented in \cite{fu2023efficient}. The proposed approach is based on the Generalized Mixed Multiscale Finite Element Method (GMsFEM) and utilizes an ILU(0) smoother. An adaptive approach in two-grid preconditioner construction is presented in \cite{yang2019two, yang2022adaptive, fu2024adaptive} for problems with heterogeneous properties.
A two-grid domain decomposition preconditioner for multiscale flows in high-contrast media is presented in \cite{galvis2010domain, efendiev2009multiscale, galvis2010domain}, where the coarse space is constructed using selected eigenvectors of a local spectral problem. 
A similar approach was proposed in the spectral algebraic multigrid method ($\rho$AMG) to solve local spectral problems in non-overlapped and overlapped local domains. It was shown that the method could give a great performance for elliptic problems with high-contrast \cite{brezina1999iterative, chartier2003spectral, chartier2007spectral}. 
The multiscale spectral multigrid solver for high-contrast flow problems with nested coarse spaces is proposed in \cite{efendiev2012multiscale, Efen_GVass_11}.
In \cite{hajibeygi2008iterative, hajibeygi2011adaptive, hajibeygi2011hierarchical},  the multiscale finite-volume method (MsFVM) is extended to an efficient iterative algorithm that converges to the fine-scale numerical solution. The main ingredients in iterative MsFVM contain recalculating basis functions with updated global boundary conditions and separate treatment of fractures with a splitting approach. In GMsFEM, we overcome these issues by proposing a systematic approach to basis construction and introducing local spectral basis functions for fractured porous media. 

In \cite{akkutlu2016multiscale, akkutlu2017multiscale}, we considered the construction of a multiscale model order reduction technique for shale gas transport in fractured porous media. We used linearization from the previous time step, which leads to updating a problem operator at each time step. For multiscale solvers with spectral coarse spaces, updating the preconditioner every time step can be very expensive. In this work, we propose a construction of an efficient linearly implicit time approximation scheme by carefully defining a linear operator to treat implicitly that is \emph{fixed for all time} and captures the stiff behavior, while treating the resulting nonlinear residual explicitly. Linear stability analysis provides a condition that needs to be satisfied to have a stable solution. 
The second purpose of the paper is to study the performance of a spectral adaptive multiscale space in a two-grid algorithm for shale gas transport in fractured heterogeneous porous media. Accurate approximation of fractured porous media on a coarse grid is challenging due to complex heterogeneities and the high-contrast nature between fracture and porous matrix permeabilities. To address this issue, we construct multiple multiscale basis functions, where each multiscale basis function represents a continuum, and there is no need for coupling terms between these continua. The basis functions for each continuum are automatically identified by solving local spectral problems. This approach showed good accuracy in constructing reduced order models for complex problems in fractured porous media in our previous works \cite{akkutlu2015multiscale, akkutlu2017multiscale}. We show that the choice of number of basis functions can be made adaptively and automatically based on a threshold of local eigenvalues. Moreover, the basis functions are associated with a coarse scale continua that can be interpreted as a coarse scale multicontinuum model \cite{chung2017non, vasilyeva2019upscaled, vasilyeva2020learning, vasilyeva2023efficient}. In this work, we show that the constructed coarse space with an appropriate choice of pointwise smoother can be used as an adaptive two-grid precondition and give a contrast-independent convergence for an iterative solver for complex shale gas transport in fractured porous media. The adaptive coarse space naturally chooses the number of coarse basis functions needed over each subdomain. Construction of such robust solvers has a relatively significant setup cost in computing basis functions, but in combining with our proposed explicit-implicit method, with linear implicit operator fixed for all time, this setup cost only needs to be done once, and can be done offline. Numerical results demonstrate robust convergence of the solver for a smaller and larger number of fractures with homogeneous and heterogeneous porous matrix properties and varying values of fracture permeability.

The paper is organized as follows. In Section 2, we describe problem formulation for shale gas transport in fractured media. In Section 3, we present a semi-discrete problem on a fine grid with finite element approximation for a mixed-dimensional formulation. Moreover, we present a detailed construction of the implicit-explicit scheme for shale gas transport processes and give a discrete formulation of the problem. In Section 4, we present a two-grid algorithm, where the coarse grid approximation is built using the adaptive spectral multiscale method and based on the spectral equivalency of the chosen smoother to diagonal part of the considered matrix. An accurate multiscale space for coarse grid approximation is a central component that gives a contrast-independent convergence for the iterative method with a two-grid preconditioner, and is related to the convergence theory of two-level iterative methods. In Section 5, we present numerical results for two test geometries with 30 and 160 fracture networks. We demonstrate robust convergence for homogeneous and heterogeneous porous matrices. The paper ends with a conclusion.

\section{Mathematical model}

To model transport processes in organic-rich shales, we incorporate the gas flow characteristics in the ultra-tight shale matrix, such as the
nonlinear gas di?usion and desorption  \cite{akkutlu2012multiscale, akkutlu2017multiscale}. Because shale formations are characterized by low porosity and extremely low permeability, the permeability enhancement by hydraulic fracturing is commonly used to create artificial pathways and extract shale gas.  
A background media (shale matrix) contains organic (kerogen) and inorganic pores. We let $c_i$ be the amount of gas in inorganic matter, and $c_k$ and $c_{\mu}$ be the amounts of free gas and adsorbed gas in kerogen per kerogen pore volume and kerogen per kerogen solid volume, respectively. An organic matrix with larger pores can store a significant amount of gas and has a viscous force. 
Gas transport within the kerogen occurs through free and adsorbed phases. The dominant transport mechanism in this nanoporous organic material is molecular diffusion, driven by concentration gradients and significantly influenced by the complex pore geometry and adsorption-desorption dynamics.

\subsection{Gas transport in shales}

Let $\phi$ be a total matrix porosity, $\varepsilon_{ks}$ be the total organic content in terms of organic grain volume per total grain volume, and $\varepsilon_{kp}$ be the kerogen pore volume per total matrix pore volume \cite{akkutlu2012multiscale, fathi2012mass, fathi2014multi, akkutlu2017multiscale}.  
For free and adsorbed gas in the kerogen, we have
\begin{equation}
\label{eq:m-k1}
\phi \varepsilon_{kp} \frac{\partial c_k}{\partial t} +
(1 - \phi) \varepsilon_{ks}  \frac{\partial c_{\mu}}{\partial t}  
- \nabla \cdot (\phi \varepsilon_{kp}  D_k \nabla c_k) 
- \nabla \cdot ((1 - \phi) \varepsilon_{ks} D_s \nabla c_{\mu}) = 0,
\end{equation}
where $D_k$ is the coefficient of diffusive molecular transport for the free gas, and $D_s$ is the coefficient of diffusive molecular transport for the adsorbed gas.
For adsorbed gas in kerogen under equilibrium conditions, we have the Langmuir isotherm
\begin{equation}
\label{eq:m-lang}
c_{\mu} = F(c_k), 
\quad 
F(c_k) = c_{\mu s} \frac{K c_k}{1 + K c_k}, 
\quad 
F'(c_k) = c_{\mu s} \frac{K}{(1 + K c_k)^2},
\end{equation}
where $c_{\mu s}$ is maximum monolayer gas adsorption on the internal kerogen solid interfaces, and $K$ is the equilibrium distribution coefficient, $K = ZRT/p_L$ ($p_L$ is the Langmuir pressure, $Z$ is the gas compressibility factor, $R$ is the universal gas constant and $T$ is the temperature).

Substituting \eqref{eq:m-lang} into  \eqref{eq:m-k1}, we obtain
\begin{equation}
\label{eq:m-k3}
(\phi \varepsilon_{kp} + (1 - \phi) \varepsilon_{ks} F'(c_k)) \frac{\partial c_k}{\partial t} -
\nabla  \cdot ( (\phi \varepsilon_{kp} D_k + (1 - \phi) \varepsilon_{ks}D_s F'(c_k) ) \nabla c_k) = 0.
\end{equation}
For the free gas mass balance in the inorganic matrix \cite{akkutlu2012multiscale}, we have
\begin{equation}
\label{eq:m-in1}
\phi (1-\varepsilon_{kp}) \frac{\partial c_i}{\partial t} -
\nabla \cdot (\phi (1-\varepsilon_{kp}) D_i \nabla c_i) -
\nabla \cdot \left(c_i \frac{\kappa_i}{\mu} \nabla p_i \right) = 0, 
\end{equation}
where $\phi_i = (1-\varepsilon_{ks}) \phi$ is the inorganic porosity, $D_i$ is the coefficient of diffusive molecular transport in the inorganic matrix, $\kappa_i$ is the inorganic matrix absolute permeability, $\mu$ is the dynamic gas viscosity, $p_i$ is the inorganic matrix pressure with $p_i = c_i ZRT$ (compressibility equation of state) \cite{fathi2009matrix}. 

We let $c_m$ be the block-average gas concentration and $c_f$ be the free gas in the fracture network.  
Then by combining \eqref{eq:m-k3} with  \eqref{eq:m-in1}, we obtain the following equation for the gas amount $c_m$ in shale matrix
\begin{equation}
\label{eq:m-m}
(\phi + (1-\phi) \varepsilon_{ks} F'(c_m)) \frac{\partial c_m}{\partial t} -
\nabla \cdot \left(   \left(\phi D + (1-\phi) \varepsilon_{ks} F'(c_m) D_s + c_m ZRT \frac{\kappa_m}{\mu} \right) \nabla c_m \right) = 0, 
\end{equation}
with $k_m = k_i$ and $D = \varepsilon_{kp} D_k + (1-\varepsilon_{kp}) D_i$. 
For free gas in a fracture network, the molecular diffusion in the fractures is small compared to the viscous flow and can be ignored \cite{akkutlu2012multiscale, fathi2012mass, fathi2014multi}. Therefore, we have the following transport problem in fractures
\begin{equation}
\label{eq:m-f1}
\phi_f \frac{\partial c_f}{\partial t} -
\div (c_f ZRT \frac{\kappa_f}{\mu} \nabla c_f) = 0.
\end{equation}
where $\phi_f$ is the fracture porosity and $\kappa_f$ is the fracture absolute permeability.
Next, we introduce the gas transport between porous matrix and fractures to write a coupled nonlinear model for multiscale gas transport in shales.

\subsection{Problem formulation}

We let $\Omega \in \mathcal{R}^d$ be a domain of porous matrix and $\gamma \in \mathcal{R}^{d-1}$ be a lower-dimensional domain corresponding to fractures. In this work, we utilize a discrete fracture model, where fractures are represented using a lower dimensional model \cite{Quarteroni2008coupling, formaggia2014reduced, martin2005modeling, vasilyeva2019upscaled}. 
Here, we have \eqref{eq:m-m} for gas transport in $\Omega$ and \eqref{eq:m-f1} for gas transport in fractures $\gamma$. By introducing coupling terms $Q_{mf} = q_{mf} (c_m - c_f)$ and $Q_{fm} = q_{fm} (c_f - c_m)$ for transport between shale matrix and fractures \cite{akkutlu2017multiscale}, we obtain the following nonlinear coupled mixed-dimensional model of shale gas transport  in heterogeneous fractured porous media: 
\begin{equation}
\begin{split}
\label{mm1}
&
a_m(x, c_m) \frac{\partial c_m}{\partial t} 
- \nabla \cdot \left(b_m(x, c_m) \nabla c_m \right) 
+ \sigma_{mf}(x, c_m, c_f) (c_m - c_f)  
= f_m(x), \quad x \in \Omega, \\
&
a_f(x) \frac{\partial c_f}{\partial t}  
- \nabla \cdot (b_f(x, c_f) \nabla c_f) 
- \sigma_{fm} (x, c_m, c_f) (c_m - c_f) 
= f_f(x), \quad x \in \gamma, \\
\end{split}
\end{equation}
with
\[
\begin{split}
&
a_m(x, c_m)
= \phi + (1-\phi) \varepsilon_{ks} F'(c_m), \quad 
a_f(x)  = \phi_f(x), \\
& 
b_m(x, c_m) 
= \phi D
+ (1-\phi) \varepsilon_{ks} F'(c_m) D_s  
+ c_m ZRT \frac{\kappa_m}{\mu}, 
\quad 
b_f(x, c_f) = c_f ZRT \frac{\kappa_f}{\mu}, 
\end{split}
\]
and 
\[
\sigma_{fm}(x, c_m, c_f) = \sigma_{mf}(x, c_m, c_f) \approx b_m(x, c_m) \zeta_{mf},
\]
where $\zeta_{mf}$ is the geometric scaling parameter representing the fracture-matrix interface and distance from matrix to fracture. 

The coupled system  of equations \eqref{mm1} is considered with the given initial condition 
\[
c_m = c_f = c_0, \quad t = 0
\]
and zero flux boundary conditions 
\[
-b_m  \nabla c_m \cdot n = 0, \quad x \in \partial \Omega, \quad 
-b_f  \nabla c_f \cdot n = 0, \quad x \in \partial \gamma,  \quad
t > 0, 
\]
where $n$ is the outer normal vector to the domain boundary.

\section{Discrete problem}

We start with an approximation by space using the finite element method. Then, we present an additive partition of the operators, and propose an implicit-explicit integration scheme. 

\subsection{Approximation by space}

With $V_m = H^1(\Omega)$ and  $V_f = H^1(\gamma)$, we have following weak form: find $(c_m, c_f) \in V_m \times V_f$ such that 
\begin{equation}
\begin{split}
\label{dm1}
 \int_{\Omega} a_m(x, c_m)& \frac{ \partial c_m}{\partial t}  v_m \, dx 
 + \int_{\Omega} b_m(x, c_m) \nabla c_m \cdot \nabla v_m \, dx \\
&  + \int_{\gamma} \sigma_{mf}(x, c_m, c_f) (c_m - c_f) v_m \, ds
 = \int_{\Omega} f_m(x)  v_m \, dx,  \quad \forall v_m \in V_m,
 \\ 
 \int_{\gamma} a_f(x, c_f) &\frac{ \partial c_f}{\partial t}  v_f \, ds 
+ \int_{\gamma}  b_f(x, c_f) \nabla c_f \cdot \nabla v_f \, ds \\
&- \int_{\gamma} \sigma_{fm}(x, c_m, c_f) (c_m - c_f) v_f \, ds 
 = \int_{\gamma} f_f(x)  v_f \, ds, \quad \forall v_f \in V_f.
\end{split}
\end{equation}
\begin{figure}[h!]
  \centering
  \includegraphics[width=0.6 \textwidth]{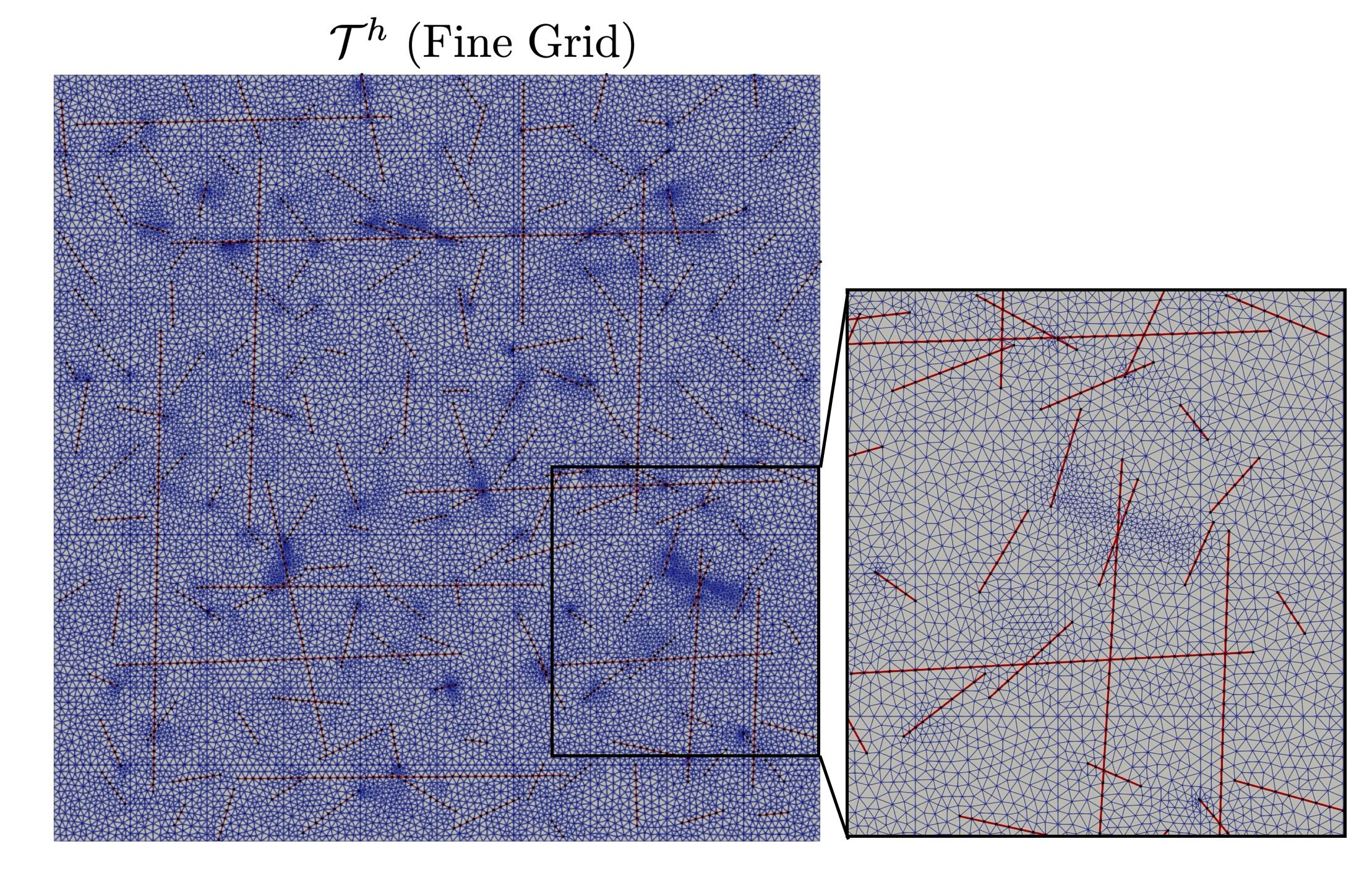}
  \caption{Illustration of a fine-scale triangulation conformed with fracture grid.}
  \label{finegrid}
\end{figure}
To construct an accurate approximation for the mixed-dimensional model, we use a discrete fracture model and construct a triangulation $\mathcal{T}^h$ of the domain $\Omega$ that is conforming with lower dimensional fractures geometries (see Figure \ref{finegrid}). 
For simplicity, we consider $\Omega \in \mathcal{R}^2$ and set $\gamma = \cup \ \mathcal{E}_{\gamma}$ as the fracture facets, where $\mathcal{E}_{\gamma} \subset \mathcal{E}_h$ and $\mathcal{E}_h$ is the all facets of the mesh $\mathcal{T}^h$.

Let $V_m^h \subset V_m$, $V_f^h \subset V_f$ and 
\[
c_m = \sum_{i=1}^{N_{\Omega}} c_{m,i} \phi_i(x), \quad 
c_f = \sum_{i=1}^{N_{\gamma}} c_{f,i} \hat{\phi}_i(x), 
\]
where $N_{\Omega}$ and $N_{\gamma}$ are the dimension of $V^h_m$ and $V_f^h$, and $\phi_i$ and $\hat{\phi}_i$ are the basis functions on fine grid. 
For $c = (c_m, c_f) \in V^h_m \times V_f^h$, we have the following coupled nonlinear system of equations
\begin{equation}
\label{dm3a}
S(c) \frac{\partial c}{ \partial t} + D(c) c = F,
\end{equation}
with
\[
S = 
\begin{pmatrix}
S_m	& 0 \\
0 &  S_f
\end{pmatrix}
, \quad 
D =
\begin{pmatrix}
L_m + Q_{mf}	& -Q_{mf} \\
-Q_{fm} & L_f + Q_{fm} 
\end{pmatrix}
, \quad 
F  =
\begin{pmatrix}
F_m \\
F_f
\end{pmatrix},
\]
where $L_m$ and $L_f$ are the stiffness matrices, and $S_m$ and $S_f$ are the mass matrices, defined as:
\begin{alignat*}{2}
S_m &= \{a^m_{ij} =  \int_{\Omega} a_m(x, c_m) \phi_i \phi_j \ dx \}, \quad 
&&S_f = \{a^f_{ij} = \int_{\gamma} a_f(x, c_f) \hat{\phi}_i \hat{\phi}_j \ ds \},
\\
L_m &= \{b^m_{ij} = \int_{\Omega} b_m(x, c_m) \nabla \phi_i \cdot \nabla \phi_j \ dx  \}, \quad 
&&L_f = \{b^f_{ij} = \int_{\gamma} b_f(x, c_f) \nabla \hat{\phi}_i \cdot \nabla \hat{\phi}_j \ ds \},  
\\
F_m &= \{f^m_j = \int_{\Omega} f_m(x) \phi_j dx \}, \quad 
&&F_f  = \{  f^f_j = \int_{\gamma} f_f(x) \hat{\phi}_j ds \}.
\end{alignat*}
For the transfer term with the same order of polynomials in continuous Galerkin finite element approximation, we have  
\[
Q_{fm} = \{q_{li} = \int_{\gamma} \sigma_{fm}(x, c_m, c_f) \hat{\phi}_i \hat{\phi}_l \ ds \},   
\]
and $Q_{mf} = Q_{fm}^T$. 
We note that we obtain a symmetric positive semi-definite matrices in approximation.

\subsection{Approximation by time}

Next, we develop a linearly stable first-order implicit-explicit scheme for approximation in time. 
In shale gas transport, our nonlinear coefficients are highly dependent on the current gas distribution in fractured media. The numerical implementation thus requires applying nonlinear iterations, such as the Newton?Raphson scheme or Picard iterations. 
A nonlinear implicit scheme would thus require updating the linearized system operator at least each time iteration and potentially each nonlinear iteration. In the context of multiscale FEMs, we would then have to reconstruct multiscale basis functions as well.
An efficient algorithm associated with problem linearization can be constructed using an additive partition of the problem operators $S$ and $B$. An additive partition allows the separation of the linear part of operators and applies an implicit scheme for a stable solution scheme. The remaining nonlinearity is moved to a previous time by approximating it explicitly from the previous known solution. With careful choice of linear component that is independent of time and nonlinear iteration, the algorithm leads to an efficient implementation that can construct a preconditioner or matrix factorization once, which is then reused through the entire simulation. This is particulaly important for multiscale finite elements, where the construction of multiscale basis functions is the dominant cost. 

For the coupled system of equations \eqref{dm3a}, we use the following additive representation of the problem operator:
\[
S(c) = S^{(lin)} + S^{(nl)}(c), \quad 
D(c) = D^{(lin)} + D^{(nl)}(c),
\]
where $S^{(lin)}$ and $D^{(lin)}$ are linear operators, that we will specify, and $S^{(nl)}$ and $D^{(nl)}$ are the nonlinear residuals:
\[
S^{(nl)}(c) = S(c) - S^{(lin)}, \quad 
D^{(nl)}(c) = D(c) - D^{(lin)}.
\]
This additive representation separates the nonlinear part by explicit approximation and allows us to work with a \emph{fixed} linear operator within a stable implicit scheme. This is in the spirit of explicit-implicit-null methods, where we add and subtract a chosen linear operator from the equation to improve stability while minimizing computational cost. 
For the coupled nonlinear system \eqref{dm3a}, we have
\[
(S^{(lin)} + S^{(nl)}(c)) \frac{\partial c}{\partial t} + (D^{(lin)} + D^{(nl)}(c)) c = F. 
\] 

Let $\tau$ be a time step size and $c^n = c(t^n)$, where $t^n = j \cdot \tau$ with $j=1,...,N_t$ and $\tau = T_{max}/N_t$. 
We construct a first-order implicit-explicit scheme and obtain the following discrete problem
\begin{equation}
\label{ein}
S^{(lin)} \frac{c^{n+1} - c^n}{\tau} 
+ S^{(nl)}(c^n) \frac{c^{n} - c^{n-1}}{\tau} 
+ D^{(lin)} c^{n+1} + D^{(nl)}(c^n) c^n = F.
\end{equation}
Due to the nonlinear interaction of the time derivative with the solution, $S(c)\partial c/\partial t$, higher order schemes as in, e.g. \cite{ascher1995implicit,ascher1997implicit,kennedy2003additive,kennedy2019higher} are not trivial and is a topic for future work. Next, we consider stability for the proposed linearized scheme \eqref{ein}.  For completeness, we include a proof of this result in the Appendix \ref{secA1}.

\begin{theorem}
\label{t:t1}
For symmetric positive definite $S^{(lin)}, S^{(nl)}(c^n),
D^{(lin)}, D^{(nl)}(c^n) > 0$, suppose the following holds:
\begin{equation}
\label{eq:stabcond}
S^{(lin)} > S^{(nl)}(c^n), \quad 
D^{(lin)} > D^{(nl)}(c^n), \quad 
\forall n=1,\ldots,N_t.
\end{equation}
Then the discrete scheme in \eqref{ein} is unconditionally stable.
\end{theorem}

Next, we consider how to accurately choose the linear part of the operator in additive representation to satisfy the stability condition from Theorem \ref{t:t1} for the shale gas transport problem. 
To have a stable time approximation scheme, we must satisfy  conditions \eqref{eq:stabcond}
with 
\[
S^{(lin)} + S^{(nl)}(c^n) = 
S^n = 
\begin{pmatrix}
S^n_m	& 0 \\
0 &  S_f
\end{pmatrix}, 
\]\[
D^{(lin)} + D^{(nl)}(c^n) =
D^n = 
\begin{pmatrix}
L^n_m + Q^n_{mf}	& -Q^n_{mf} \\
-Q^n_{fm} & L^n_f + ^nQ_{fm} 
\end{pmatrix}, 
\]
where
\[
\begin{split}
&S^n_m = \left\{\int_{\Omega} a_m(x, c^n_m) \phi_i \phi_j \ dx \right\}, \quad 
S_f = \left\{\int_{\gamma} a_f(x) \hat{\phi}_i \hat{\phi}_j \ ds \right\},\\
&Q^n_{fm} = \left\{\int_{\gamma} \sigma_{fm}(x, c^n_m, c^n_f) \hat{\phi}_i \hat{\phi}_l \ ds \right\},  \\
&L^n_m = \left\{\int_{\Omega} b_m(x, c^n_m) \nabla \phi_i \cdot \nabla \phi_j \ dx  \right\}, \quad 
L^n_f = \left\{\int_{\gamma} b_f(x, c^n_f) \nabla \hat{\phi}_i \cdot \nabla \hat{\phi}_j \ ds \right\}.
\end{split}
\]
In the shale gas transport, we have
\[
\begin{split}
&
a_m(x, c^n_m)
= \phi + (1-\phi) \varepsilon_{ks} F'(c^n_m), \quad 
a_f(x)  = \phi_f(x),\\
&
\sigma_{fm}(x, c^n_m, c^n_f) = b_m(x, c^n_m) \zeta_{mf},
\\
& 
b_m(x, c^n_m) 
= \phi D
+ (1-\phi) \varepsilon_{ks} F'(c^n_m) D_s  
+ c^n_m ZRT \frac{\kappa_m}{\mu}, 
\quad 
b_f(x, c^n_f) = c^n_f ZRT \frac{k_f}{\mu},
\end{split}
\]
with $F'(c) = c_{\mu s} \frac{K}{(1 + K c)^2}$.

\begin{figure}[h!]
\centering
\includegraphics[width=0.48 \textwidth]{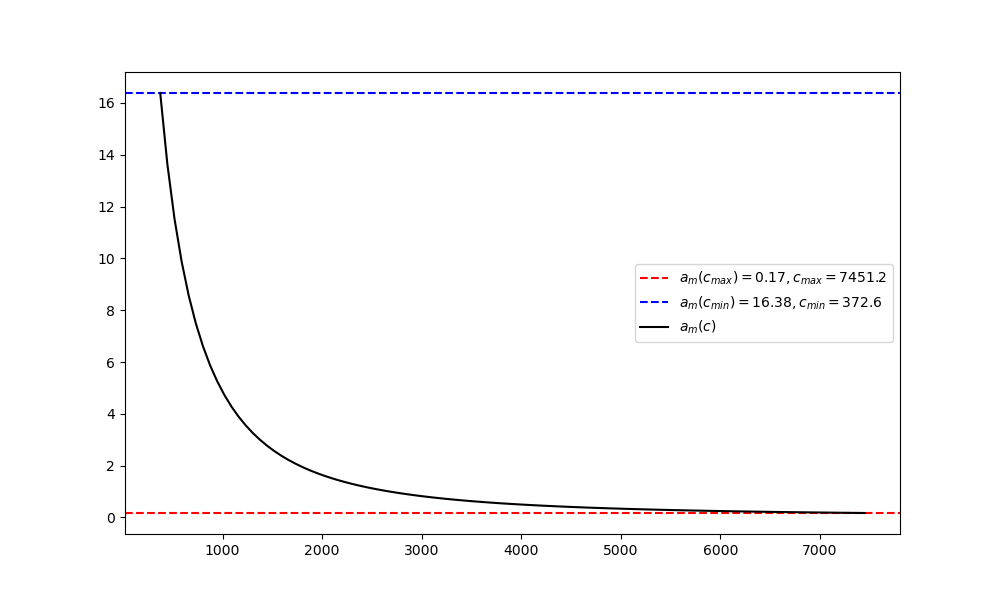} \ \
\includegraphics[width=0.48 \textwidth]{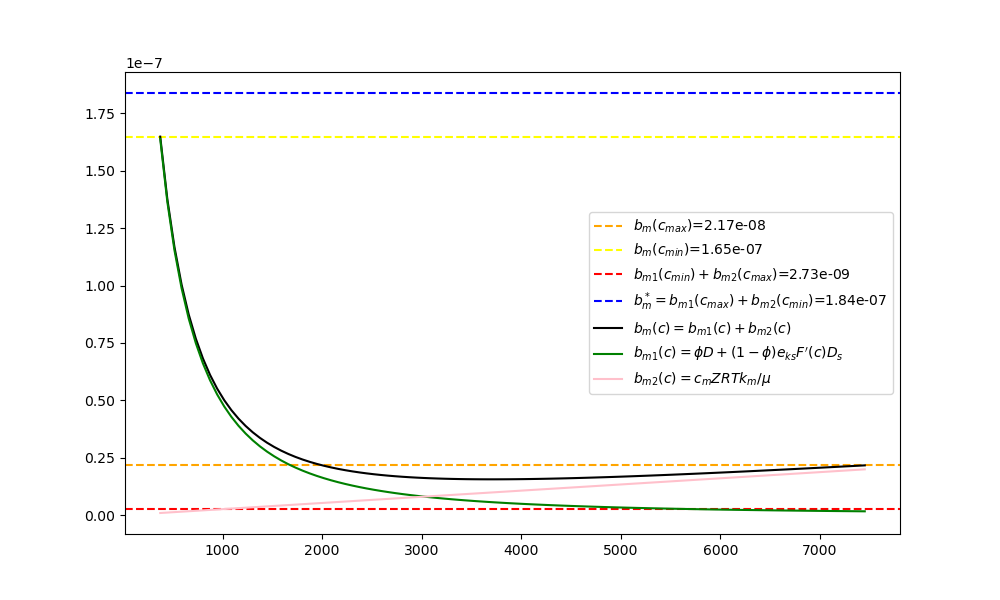}
\caption{Illustration of nonlinear coefficients $a_m(c)$ and $b_m(c)$ for shale gas transport.}
\label{fig:split}
\end{figure}

To satisfy the stability estimate, we should use an upper bound for operators $S^n$ and $D^n$ in an additive representation. For $S^n$, we have a block diagonal matrix with linear block $S_f$ and linearized block $S^n_m$ that depends on coefficients $a_f$ and $a_m$, respectively. Similarly, for operator $D$, we have a block structure with blocks that rely on coefficients $b_f$ and $b_m$. Therefore, to find a stable splitting, we should see an upper bound for coefficients
\[
a_m(x, c^n_m) \leq a_m^*(x), \quad 
b_m(x, c^n_m) \leq b_m^*(x),\quad 
b_f(x, c^n_f) \leq b_f^*(x).
\]
Let $c^n_m, c^n_f \in [c_{min}, c_{max}]$, where $c_{min}$ and $c_{max}$ are the minimum and maximum gas amounts. For example, the minimum is given at the boundary or as concentration (pressure) on the production well, and the maximum is an initial concentration. 
Therefore, for a given nonlinearity of $a_m(x, c^n_m) =  \phi(x) + (1-\phi(x)) \varepsilon_{ks} F'(c^n_m)$, we obtain maximum value at minimum concentration, $a_m^*(x) = \phi(x) + (1-\phi(x)) \varepsilon_{ks} (F')^* = a_m(x, c_{min})$  and $(F')^* = F'(c_{min})$ (see left Figure \ref{fig:split}). 
A nonlinear coefficient $b_m$, we have composite nonlinear behavior that depends on absorbed gas in kerogen and a second part that depends on gas flow, $b_m(x, c^n_m) = b_{m1}(x, c^n_m) + b_{m2}(x, c^n_m)$. Similarly to $a_m$, for $b_{m1}(x, c^n_m) = \phi D + (1-\phi) \varepsilon_{ks} F'(c^n_m) D_s$, we have maximum value at minimum concentration, $b_{m1}^*(x) = \phi D + (1-\phi) \varepsilon_{ks} (F')^* = b_{m1}(x, c_{min})$  (see Figure \ref{fig:split}, green line on right plot). For $b_{m2}(x, c^n_m) = c^n_m ZRT \frac{\kappa_m}{\mu}$, we have maximum value at maximum concentration, $b_{m2}^*(x) = b_{m2}(x, c_{max})$ (see Figure \ref{fig:split}, pink line on right plot). Similarly, we have $b_f^*(x) = b_f(x, c_{max})$ for fractures. 
Finally, we set
\[
a_m^*(x) = a_m(x, c_{min}), \quad 
a_f^*(x) = a_f(x), 
\]\[
b_m^*(x) = b_{m1}(x, c_{min}) + b_{m2}(x, c_{max}), \quad 
b_f^*(x) = b_f(x, c_{max}).
\]

Finally, we solve a linear system of equations on each time iteration. 
The discrete system of equations \eqref{ein} can be written in the following form for $c^{n+1} = (c^{n+1}_m, c^{n+1}_f)$:
\begin{equation}
\label{s1}
A c^{n+1} = b^{n}, 
\end{equation}
with
\[
\begin{split}
&A = S^{(lin)} + \tau  D^{(lin)}, \\
&b^{n} 
= S^{(lin)} c^{n} + \tau F 
- S^{(nl)}(c^n) (c^{n} - c^{n-1})
- \tau D^{(nl)}(c^n) c^n,
\end{split}
\]
where $A$ is the discrete transition matrix from time $t_{n}$ to $t_{n+1}$.

\section{Adaptive two-grid preconditioning by spectral coarse space}

Here we seek an iterative solver for the system of linear equations \eqref{s1} arising in fractured porous media with high-contrast properties, $k_f/k_m >> 1$. In particular, our goal is to design a two-grid preconditioner that converges independently of the coefficient contrast $k_f/k_m$. Such contrast can be challenging for geometric multigrid methods, where coarse-grid problems that are unable to capture the effects of sharp contrast or heterogeneity can result in a poor a coarse-grid correction.

In order to construct an efficient preconditioner, we appeal to pointwise smoothing coupled with a spectral coarse space. Spectral coarse spaces have seen remarkable success in two-level methods, where transfer operators to and from a Galerkin coarse grid are defined based on some local (generalized) eigenvalue problem, e.g. \cite{egvdd20,daas2024robust,al2023efficient,barker2017spectral}. Here we avoid the more robust domain-decomposition/(overlapping) Schwartz approach for fine-grid smoothing used in, e.g. \cite{daas2024robust,al2023efficient} to reduce the computational cost, and because we demonstrate we can obtain good results with a cheap pointwise approach. For constructing a spectral coarse space, there are different approaches in the literature. Here we construct a coarse space as an adaptive spectral multiscale space for fractured porous media in the spirit of GMFEMs, which have shown success for high-contrast fractured porious media \cite{akkutlu2015multiscale, akkutlu2018multiscale}. In contrast with classic geometric or algebraic restriction and interpolation operators, our transfer operators are built on local spectral problems, and lead to a very accurate multicontiuum coarse approximation \cite{akkutlu2015multiscale, akkutlu2017multiscale, vasilyeva2024generalized} with accurate adaptive fracture approximation.

\subsection{Two-grid methods}

Following equation \eqref{s1}, for each implicit time step we solve a system of linear equations
\[
A y = b,
\]
with
\begin{equation}\label{eq:mat}
A = 
\begin{pmatrix}
A_{mm} & A_{mf} \\
A_{fm} & A_{ff}
\end{pmatrix}, \quad 
y = 
\begin{pmatrix}
y_m \\
y_f 
\end{pmatrix}, \quad 
b = 
\begin{pmatrix}
b_m \\
b_f 
\end{pmatrix}.
\end{equation}
Here, the full matrix $A = A^T \geq 0$ is symmetric positive semi-definite, where $A \in \mathcal{R}^{N_h \times N_h}$, $y \in \mathcal{R}^{N_h}$, and $b \in \mathcal{R}^{N_h}$, with $N_h = N_{\Omega} + N_{\gamma}$, for total fine-grid degrees of freedom for porous matrix, $N_{\Omega}$, and fractures, $N_{\gamma}$.

Our two-grid method follows the framework of Algebraic Multigrid Method (AMG) methods \cite{ruge1987algebraic, xu2017algebraic}, consisting of three steps, (i) a pointwise stationary linear iterative method (pre-smoothing) on the fine-grid system, (ii) coarse-grid correction via projecting the residual to be approximated by some coarse operator, and (iii) a similar post-smoothing iteration on the fine grid. The efficiency of the approach depends on a proper complementary relationship between the smoother and the coarse-grid correction \cite{vassilevski2008multilevel, falgout2005two, falgout2004generalizing, brezina2011smoothed, vassilevski2011coarse}, which the high-contrast nature of fractured porous media makes more challenging.

Let  $M  \in \mathcal{R}^{N_h \times N_h}$ be a smoothing operator and $A_H = P^T A P \in \mathcal{R}^{N_H \times N_H}$ be a coarse grid matrix constructed using interpolation and restriction operators, $P \in \mathcal{R}^{N_H \times N_h}$ and $R = P^T \in \mathcal{R}^{N_h \times N_H}$. 
Then, for a given initial guess $y^{(0)}$, we define the symmetric two-grid algorithm as follows:
\begin{enumerate}
\item \textit{Pre-smoothing}: $y^{(1)} =  y^{(0)} + M^{-1} (b - A y^{(0)})$.
\item \textit{Coarse-grid correction}:
\begin{enumerate}
\item  \textit{Restriction}:  $r_H = P^T (b - A y^{(1)})$.
\item \textit{Coarse-grid solution}: $A_H e_H = r_H$ with $A_H = P^T A P$
\item \textit{Interpolation and update}: $y^{(2)} = y^{(1)} + P e_H$.
\end{enumerate}
\item \textit{Post-smoothing}: $y_{TG} =  y^{(2)} + M^{-T} (b - A y^{(2)})$.
\end{enumerate}
Altogether, we have the following error transfer operator of two-level AMG \cite{falgout2005two,notay2007convergence}:
\[
E_{TG} = (I - M^{-T} A) (I - P A_H^{-1} P^T A) (I - M^{-1} A).
\]
With $\nu \geq 0$ pre- and post-smoothing iterations, we obtain error propagation
\[
E_{TG}  = (I - M^{-T} A)^{\nu} (I - P A_H^{-1} P^T A) (I - M^{-1} A)^{\nu}
\]
with $\nu$ is the number of smoothing iterations. 

Smoothing iterations are typically designed to remove high-frequency errors \cite{vassilevski2008multilevel, falgout2005two, fu2023efficient, yang2019two, xu2022convergence}. Algebraically, pointwise relaxation attenuates error associated with the largest eigenvalues. A coarse grid correction is designed to attenuate the remaining error not effectively reduced by smoothing \cite{akkutlu2015multiscale, fu2023efficient, yang2019two}. In simpler and/or smoother problems, such as elliptic problems with high regularity, there is often good correspondence between algebraic smoothness (small eigenvalues) and geometric smoothness. However, high-contrast problems complicate this relationship, often requiring specialized methods in one or both of these components. In particular, algebraically smooth error (i.e. error associated with small eigenvalues) does not necessarily appear geometrically smooth due to the high-contrast fracture space. 

In this work, we use a pointwise symmetric Gauss-Seidel iteration as a pre- and post-smoother, to keep the cost of fine-grid smoothing minimal (avoiding, e.g., large overlapping block/domain solves). For $M = L$, where $L=U^T$ is the lower triangular part of the matrix, this method produces an SPD symmetric smoother denoted as $\bar{M} = M (M + M^T  - A)^{-1} M^T$, which is spectrally equivalent to the diagonal of matrix $A$ \cite{vassilevski2008multilevel}. Thus assuming smoothing to attenuate error associated with large eigenvalues, we rely on coarse-grid correction to represent error associated with small eigenvalues, although such error is not always geometrically smooth.

\subsection{Adaptive spectral space} 

\begin{figure}[!hbt]
  \centering
  \includegraphics[width=0.6 \textwidth]{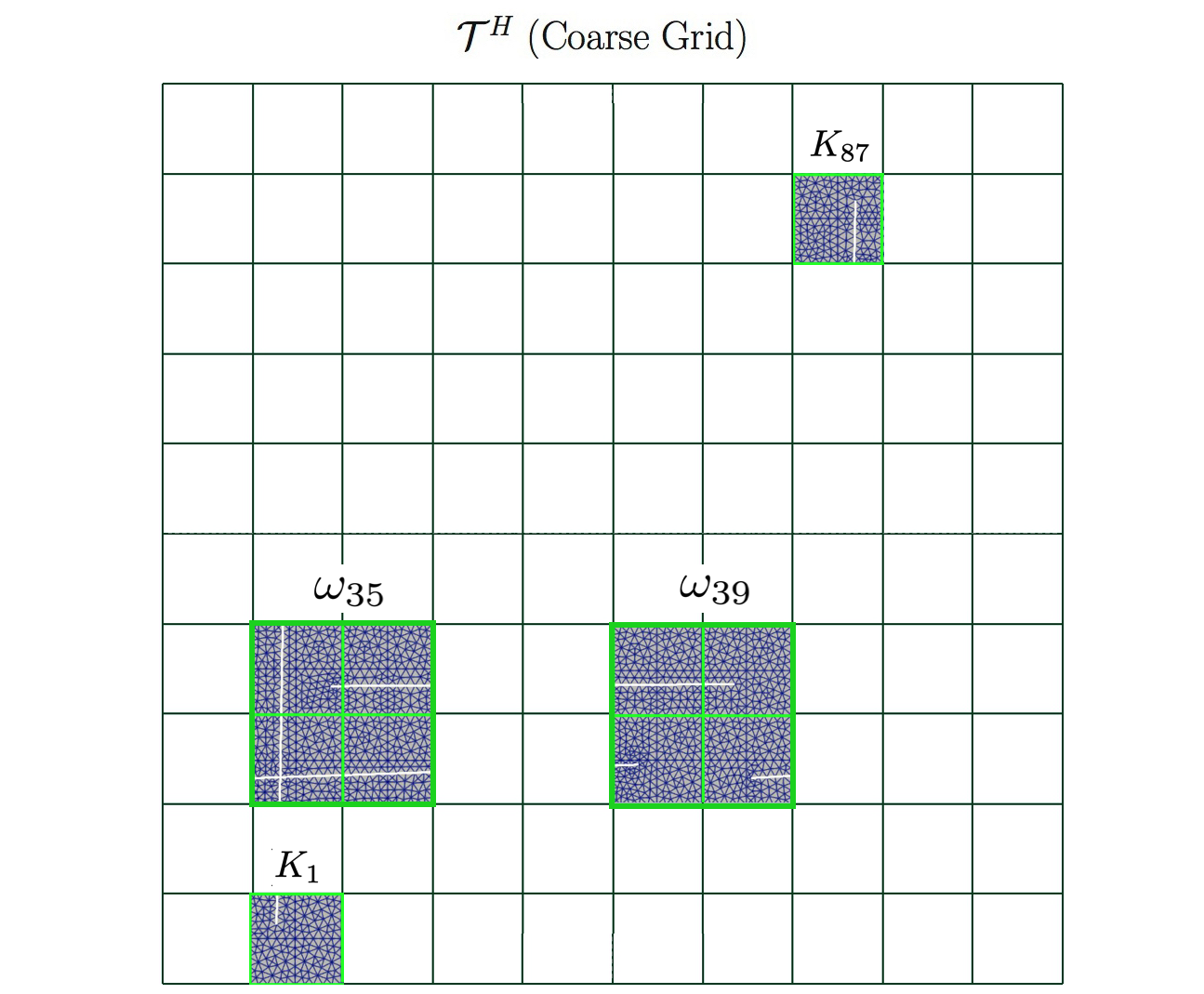}
  \caption{Illustration of a coarse grid $\mathcal{T}^H$ with local domains ($\omega_{35}$ and $\omega_{39}$) and coarse elements ($K_1$ and $K_{87}$).}
  \label{fig:ms-coarse}
\end{figure}

We consider the following notation to define fine and coarse problems:
\begin{enumerate}
\item \textit{Fine grid}: let $\mathcal{T}^h$ denote the fine-level triangulation with discrete fracture approximation using finite element method, $V_h = V^h_m \times V^h_h$ and $N_h = dim(V_h)$.
\item \textit{Coarse grid}: let $\mathcal{T}^H$ denote the coarse grid with spectral multiscale space $V_H$ and adaptive choice of basis functions, $N_H = dim(V_H)$.
\end{enumerate}
To construct an accurate approximation for the mixed-dimensional model on the coarse grid, we construct an adaptive multiscale space based on the Generalized Multiscale Finite Element Method (GMsFEM) \cite{efendiev2013generalized, akkutlu2015multiscale}. 
On a coarse grid $\mathcal{T}^H$, we use $\{x_i\}_{i=1}^{N}$ (where $N$ denotes the number of coarse nodes) to denote the vertices of the coarse mesh $\mathcal{T}^H$, and define the neighborhood of the node $x_i$ by $\omega_i$, where
\[
\omega_i=\bigcup\{ K_j\in\mathcal{T}^H; \ x_i\in \overline{K}_j\}.
\]
See Figure \ref{fig:ms-coarse} for an illustration of coarse neighborhoods (local domain, $\omega_i$) and coarse elements ($K_i$) related to the coarse grid.

We consider a continuous Galerkin formulation of the coarse grid approximation and define the continuous Galerkin spectral multiscale space as
\[
V_{H}  = \text{span} \{ \chi_i \psi_k^{\omega_i} : \ 1 \leq i \leq N \ \text{and} \ 1 \leq k \leq m_i  \},
\]
where $\{\chi_i\}$ are linear partition of unity functions, $\{\psi_k^{\omega_i}\}_{i=1}^{m_i}$ are the multiscale basis functions defined in local domain $\omega_i$, and subscript $k$ represents the numbering of the basis functions.
On the multiscale space, the solution can be represented as follows:
\[
u_{\text{ms}}(x, t)=\sum_{i = 1}^{N}\sum_{k=1}^{m_i} u_{k}^i(t)\phi_k^{\omega_i}(x),
\]
where $\phi_k^{\omega_i} = \chi_i \psi_k^{\omega_i}$.

In order to construct the multiscale space $V_{\text{off}}^\omega$, we solve local spectral problems in each $\omega$. The analysis in \cite{egw10, vassilevski2008multilevel, efendiev2013generalized} motivates the following generalized eigenvalue problem: 
\begin{equation} 
\label{sp}
A^{\omega_i} \psi_k^{\omega_i} = \lambda_k^{\omega_i}  D^{\omega_i}\psi_k^{\omega_i},
\end{equation}
where
\[
A^{\omega_i} = 
\begin{pmatrix}
A^{\omega_i}_{mm} & A^{\omega_i}_{mf} \\
A^{\omega_i}_{fm} & A^{\omega_i}_{ff}
\end{pmatrix}, 
\]
denotes a principle submatrix of $A$ \eqref{eq:mat} over subdomain $\omega_i$ with a modified diagonal part to preserve zero flux boundary conditions on local boundaries in continuous Galerkin finite element approximation, and $D^{\omega_i}$ is the diagonal scaling matrix given by the diagonal of $A^{\omega_i}$. In AMG theory, two-level convergence is typically analyzed with respect to the symmetrized smoother, but as discussed in \cite{vassilevski2008multilevel}, for pointwise smoothing such as symmetric Gauss-Seidel, the symmetrized smoother is spectrally equivalent to the diagonal of the matrix, so we use $D$ in \eqref{sp} for simplicity.

\begin{figure}[htb]
\centering
\begin{subfigure}{0.32\textwidth}
\centering
\includegraphics[width=1 \textwidth]{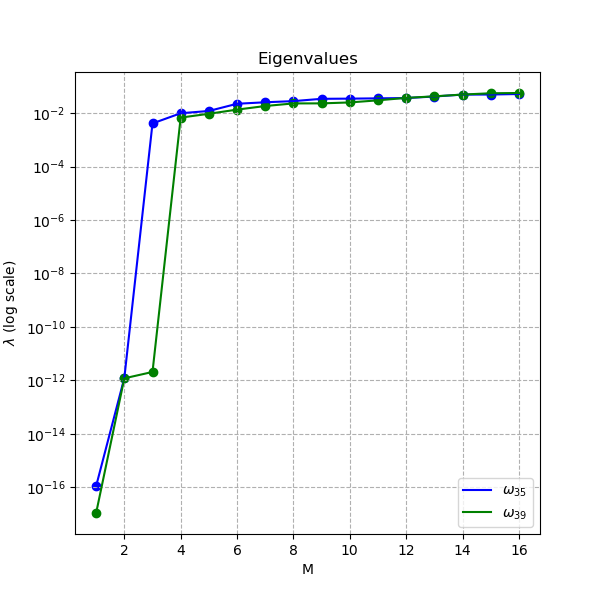}
\caption{Eigenvalues in local domains  $\omega_{35}$ and $\omega_{39}$.}
\end{subfigure}
\begin{subfigure}{0.65\textwidth}
\begin{subfigure}{1\textwidth}
\includegraphics[width=1 \textwidth]{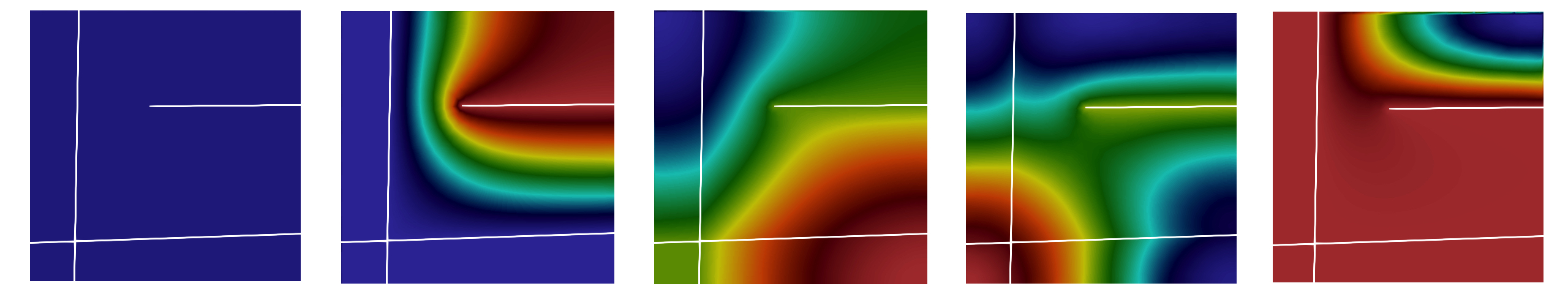}
\caption{Eigenvectors in $\omega_{35}$ with $(\lambda^{\omega_{35}}_1, \lambda^{\omega_{35}}_2, \lambda^{\omega_{35}}_3, \lambda^{\omega_{35}}_4, \lambda^{\omega_{35}}_5) = (1.12 \cdot 10^{-16},  1.23 \cdot 10^{-12},  4.08 \cdot 10^{-3},  9.83 \cdot 10^{-3}, 1.19 \cdot 10^{-2})$.}
\end{subfigure}
\begin{subfigure}{1\textwidth}
\includegraphics[width=1 \textwidth]{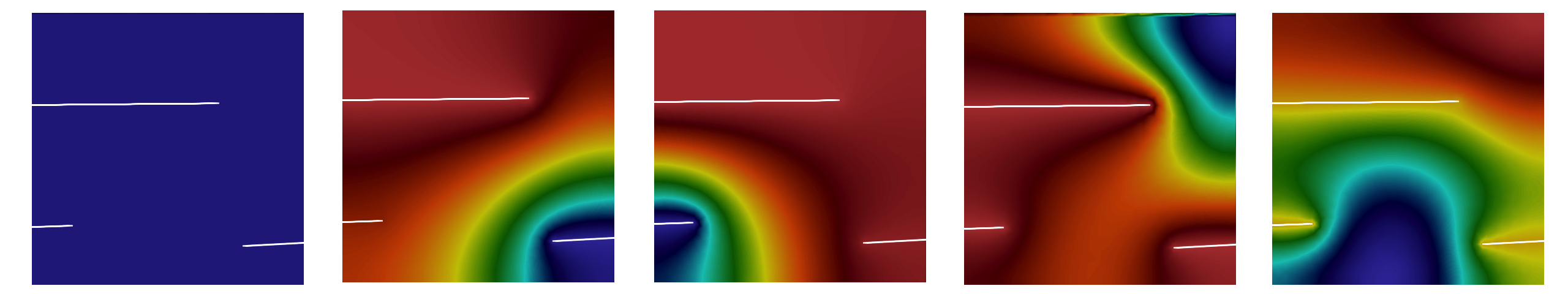}
\caption{Eigenvectors in $\omega_{39}$ with $(\lambda^{\omega_{39}}_1, \lambda^{\omega_{39}}_2, \lambda^{\omega_{39}}_3, \lambda^{\omega_{39}}_4, \lambda^{\omega_{39}}_5) = (1.10 \cdot 10^{-17},  1.20 \cdot 10^{-12},  2.07 \cdot 10^{-12},  6.75 \cdot 10^{-3},  9.39 \cdot 10^{-3})$.}
\end{subfigure}
\end{subfigure}
\caption{Eigenvalues and eigenvectors in local domains $\omega_{35}$ and $\omega_{39}$}
\label{ms}
\end{figure}

To generate the offline space, we sort eigenvalues in ascending order 
\[
\lambda_1^{\omega_i} \leq \lambda_2^{\omega_i} \leq \ldots \leq \lambda_{m_i}^{\omega_i} \leq \ldots,
\]
and choose the generalized eigenvectors associated with the smallest $m_i$ eigenvalues from equation \eqref{sp}. These represent local approximations of the $m_i$ modes least attenuated by our pointwise relaxation scheme.
In Figure \ref{ms}, we plot an illustration of the first five eigenvalues and corresponding eigenvectors for two local domains. From the eigenvalue plot, we observe the first two eigenvalues are very small, the first corresponding to the nullspace associated with natural boundary conditions, and starting from the third, we obtain the same order of eigenvalues with slow increase; therefore, we should take two basis functions in $\omega_{35}$ ($m_{35} = 2$). For the local domain $\omega_{39}$, according to the eigenvalue plot, we take three basis functions and set $m_{39} = 3$.  In this work, we choose $m_i$ based on a threshold, $m_i =\max i$ such that $\lambda_i < \delta_{\lambda}$, for some specified $\delta_\lambda$. For example, for the domains represented in Figure \ref{ms}, we choose $\delta_{\lambda} = 10^{-3}$. Like any adaptive algorithm, the choice of $\delta_{\lambda}$ remains a parameter to tune, but a single choice facilitates adaptive behavior across the global problem, as demonstrated later.

Note that generalized eigenvectors $\{\psi^{\omega_i}_k\}$ form a locally orthonormal basis with respect to inner product $v^T D^{\omega_i} u$. 
Then we define the local projection as follows 
\begin{equation}
\label{pl}
P^{\omega_i} v 
= \sum_{k=1}^{m_i}c_{i,k} {\psi}_k^{\omega_i},\quad
c_{i,k} =  (  ({\psi}_k^{\omega_i})^T D^{\omega_i}v), 
\quad v \in V_h,
\end{equation}
and we have the following inequality \cite{abreu2019convergence, efendiev2011multiscale, vasilyeva2024generalized}
\begin{equation}
\label{plest}
||v - P^{\omega_i} v ||^2_{D^{\omega_i} } 
\leq \frac{1}{ \lambda_{m_i+1}^{\omega_i} } ||v - P^{\omega_i} v ||^2_{A^{\omega_i}} 
\leq \frac{1}{ \lambda_{m_i+1}^{\omega_i} } ||v||^2_{A^{\omega_i}},
\end{equation}
with $||u||^2_{D} = u^T D u = (u, u)_{D}$ and $||u||^2_{A} = u^T A u =  (u, u)_{A}$. 

We multiply the selected eigenfunctions by partition of unity functions to create the interpolation and restriction matrices
\begin{equation}
\label{restr}
P = \left[ 
\chi_1 \psi_{1}^{\omega_1}, \ldots, \chi_1 \psi_{m_1}^{\omega_1},\ldots ,
\chi_N \psi_{1}^{\omega_N}, \ldots, \chi_N \psi_{m_N}^{\omega_N} \right], \quad 
R = P^T, 
\end{equation}
where $m_i$ denotes the number of offline eigenvectors that are chosen for each coarse node $i$. 
Note that, in forming the interpolation matrix $P$, we do a map of degrees of freedom from local domain $\omega_i$ to global domain $\Omega$, when approximation spaces for the local spectral problem in equation \eqref{sp} and the global fine-grid problem in equation \eqref{ein} are the same. 
We note that the construction yields continuous basis functions due to the multiplication of local eigenvectors with the initial (continuous) partition of unity. 

The proposed adaptive spectral space is constructed in the offline stage, and is designed to choose a minimal number of basis functions over each local subdomain. The local problems can be solved independently in each subdomain, leading to fully parallel computations in small local domains. Furthermore, in general, the choice of the coarse grid is connected with the local domain size, where we should balance fast eigenvalue solutions in the local domains and the small size of the coarse-scale system that is relatively cheap to solve directly. An extension to the multilevel setting will be considered in future work. Note that  adaptivity in coarse space construction can also be done as an online enrichment process by a-posteriori error estimate and based on the local residual in $\omega_i$ \cite{Chung_adaptive14}. The a-posteriori error indicator gives an estimate of the local error in order to add more basis functions accordingly to improve the solution. However, this approach will require the reconstruction of a interpolation operator $P$. Online adaptivity \cite{chung2015residual} can also be adopted from residual based multiscale space enrichment. Such approaches will be considered in future works for more complex nonlinear shale gas transport problems, where transfer operators will need to updated dynamically as the problem evolves.

\subsection{Convergence of two-grid algorithm}

Following \cite{vassilevski2008multilevel,falgout2005two}, two-grid error propagation can be expressed as
\[
E_{TG} = I - B^{-1}_{TG} A \quad \text{with} \quad
B^{-1}_{TG} = \bar{M}^{-1} + (I - M^{-T} A) R^T A^{-1}_H R (I - A M^{-1}).
\] 
In order to establish convergence properties of the two-grid algorithm, sufficient conditions are proving \cite{vassilevski2008multilevel, brezina2011smoothed, falgout2005two}
\[
 v^T A v \leq v^T B_{TG} v \leq K_{TG} \ v^T A v,  \quad \textnormal{where }
 K_{TG} \coloneqq \text{Cond}(B_{TG}^{-1} A),
\]
which is equivalent to 
\[
0 \leq v^T A E_{TG} v \leq \rho_{TG} v^T A v, \quad \textnormal{where }
\rho_{TG} = 1 - \frac{1}{K_{TG}}.
\]

For the Gauss-Seidel smoother, we have $A = D - N - N^T$, where $D$ is the diagonal of $A$ and $-N$ is the strictly lower triangular part of $A$. Then $M = D-N$ and 
\[
M^T + M - A = D, \quad 
\tilde{M} = (D-N^T) D^{-1} (D-N),
\]
where $\tilde{M}$ is the symmetric matrix that is spectrally equivalent to $D$ \cite{vassilevski2008multilevel}. 

Next, we write the interpolation operator $P$ as follows
\begin{equation}
\label{pg}
P v = \sum_{i=1}^{N}  \chi_i (P^{\omega_i} v),
\end{equation}
and  $v - P v  = \sum_{i=1}^N \chi_i (v - P^{\omega_i}  v)$. \

Then using properties of partition of unity functions $\chi_i$ ($\chi_i \leq 1$ and $|\nabla  \chi_i| \leq 1/H^2$), we can obtain the following estimate for the global projection \cite{melenk1996partition, abreu2019convergence, efendiev2011multiscale, vasilyeva2024generalized}
\[
||u - P u ||_{D}^2 
= \sum_K ||u - P u ||_{D^K}^2 
 \preceq
\sum_{K}  \sum_{y_i \in K}  ||\chi_i (u - P^{\omega_i}  u)||^2_{D^{K}}
 \preceq
\sum_K  \sum_{y_i \in K} ||u - P^{\omega_i}  u||_{D^{\omega_i}}^2
\]
We let  
$\lambda_{m+1} = \min_K  \lambda_{K, m+1}$ and 
$\lambda_{K, m+1} = \min_{y_i \in K}  \lambda_{m_i+1}^{\omega_i}$. 
If we map the local generalized eigenvalue problem \eqref{sp} to size one domain \cite{abreu2019convergence}, then the eigenvalues scale with $H^{-2}$ ($\Lambda^* \approx \lambda_{m+1} H^{-2}$) and combined with estimates for local projection \eqref{plest}, we obtain 
\begin{equation}
\label{pgest}
||v - P v ||_{D}^2 \leq 
\frac{H^2}{\Lambda^*}  ||v||^2_{A}.
\end{equation}

Finally, we combine the estimate \eqref{pgest} with  spectral equivalency of smoother $M$ to $D$ and obtain
\[
||v - P v ||_{\tilde{M}}^2 
\leq 
\frac{c_mH^2}{\Lambda^*} v^T A v,
\]
and we have $K_{TG} = \frac{c_mH^2}{\Lambda^*}$, where $c_m$ depends on spectral equivalence of $D$ and $\tilde{M}$.

\section{Numerical results}

We consider the model problem in fractured porous media in the domain $\Omega = [0,1]^2$ with 30 and 160 fracture lines for Test 1 and Test 2, respectively. 
We set the source term in fractures located in the lower left area of the domain ($x \in [0.1, 0.15] \times [0.05, 0.1]$) and upper right area ($x \in [0.6, 0.65] \times [0.9, 0.95]$).  

\begin{figure}[h!]
\centering
\begin{subfigure}{0.24\textwidth}
\includegraphics[width=1 \textwidth]{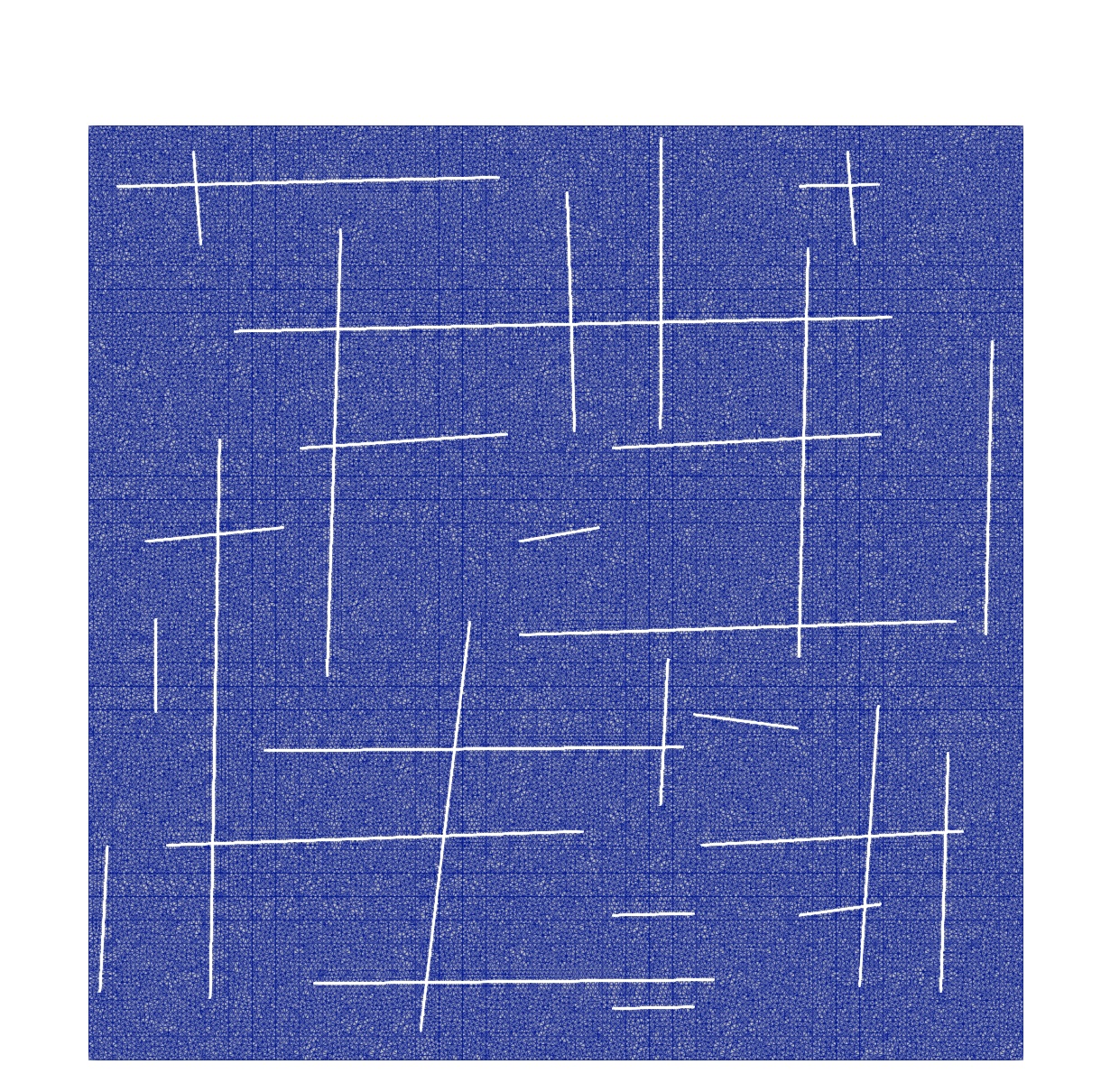}
\caption{Test 1}
\end{subfigure}
\begin{subfigure}{0.24\textwidth}
\includegraphics[width=1 \textwidth]{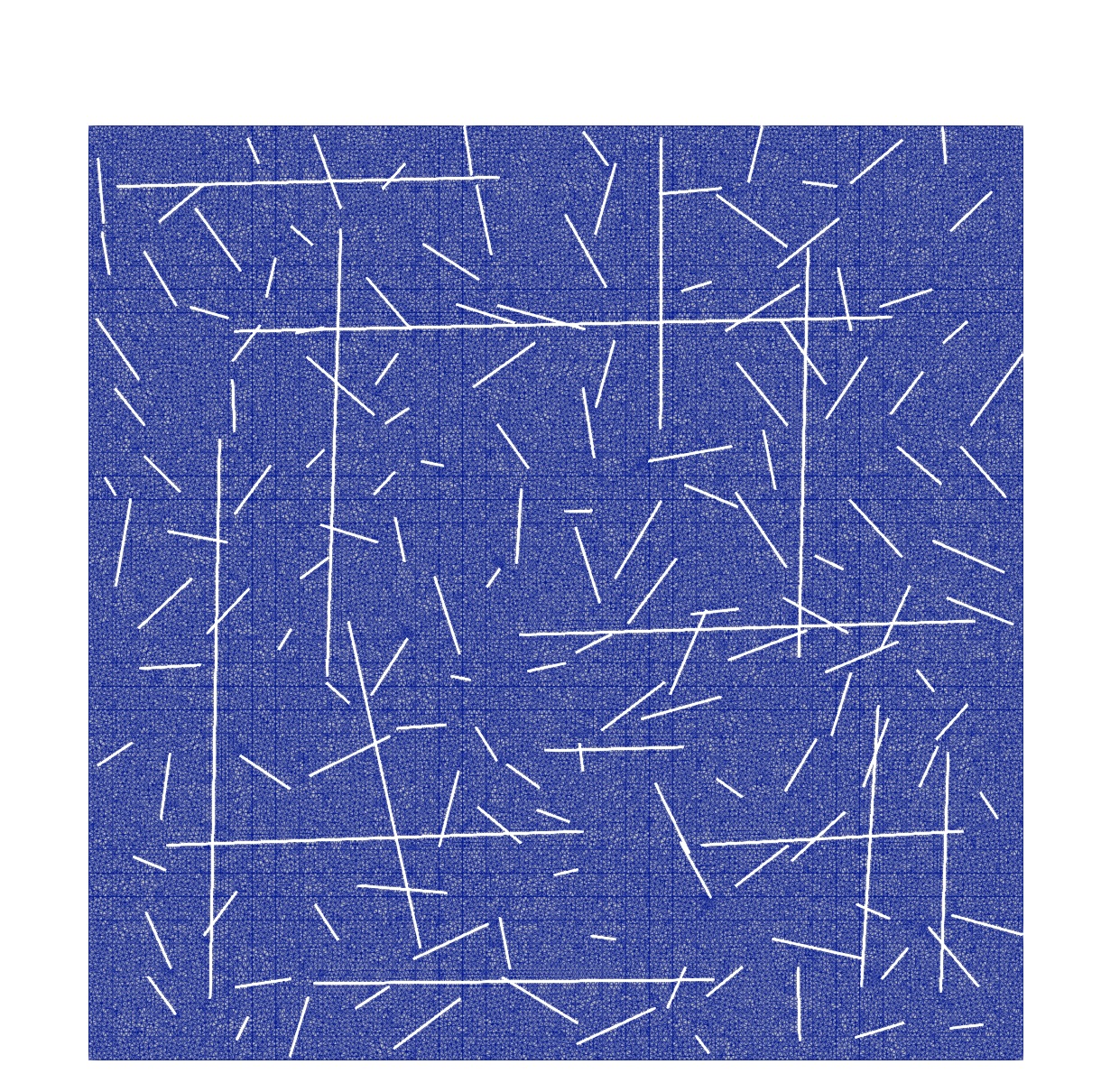}
\caption{Test 2}
\end{subfigure}
\begin{subfigure}{0.24\textwidth}
\includegraphics[width=1 \textwidth]{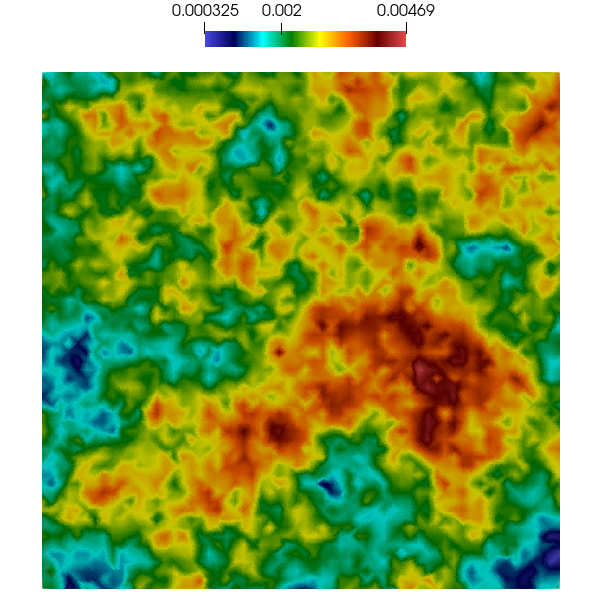}
\caption{$\phi_m$}
\end{subfigure}
\begin{subfigure}{0.24\textwidth}
\includegraphics[width=1 \textwidth]{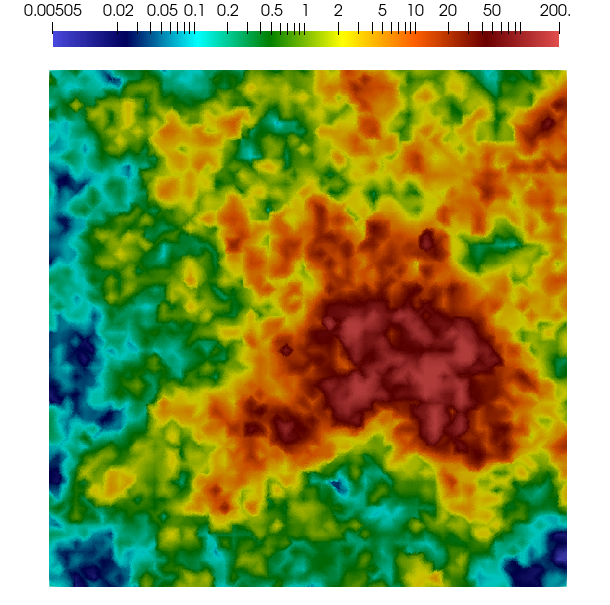}
\caption{$k_m$ (log scale)}
\end{subfigure}
\caption{Computational domains and heterogeneous coefficients. Fracture lines (white color) and fine-scale triangulation (blue color). 
First plot: Test 1 with 30 fracture lines, $N_h = 62975 + 1796$. 
Second plot: Test 2 with 160 fracture lines, $N_h = 64629 + 3481$. 
Third and fourth plots: Heterogeneous coefficients $\phi_m$ and $k_m$ for porous matrix for Test 1b and Test 2b.}
\label{fig:geom}
\end{figure}

We consider two test geometries (Figure \ref{fig:geom}): 
\begin{itemize}
\item \textit{Test 1} with 30 fracture lines. 
The fine-scale triangulation for porous matrix domain $\Omega$ contains 62975 vertices and 125148 cells (triangles). 
Fine-scale grid for fractures $\gamma$ contains 1796 vertices and 1784 cells (lines). 
\begin{itemize}
\item \textit{Test 1a}: 30 fractures with homogeneous background ($\phi = 0.02$ and $\kappa_m = 10^{-20}$ [m$^2$]).
\item \textit{Test 1b}: 30 fractures with heterogeneous background ($\phi = \phi_m \cdot 0.02$ and $\kappa_m = k_m \cdot 10^{-20}$ are depicted in Figure \ref{fig:geom}).
\end{itemize}
\item \textit{Test 2} with 160 fracture lines. 
The fine-scale triangulation for porous matrix domain $\Omega$ contains 64629 vertices and 128452 cells (triangles). 
Fine-scale grid for fractures $\gamma$ contains 3481 vertices and 3362 cells (lines). 
\begin{itemize}
\item \textit{Test 2a}: 160 fractures with homogeneous background ($\phi = 0.02$ and $\kappa_m = 10^{-20}$ [m$^2$]
\item \textit{Test 2b}: 160 fractures with heterogeneous background ($\phi = \phi_m \cdot 0.02$ and $\kappa_m = k_m \cdot 10^{-20}$ are depicted in Figure \ref{fig:geom}).
\end{itemize}
\end{itemize}
Note that the fine-scale grid for the porous matrix conforms with the fine-scale grid for fractures. For fractures, we set $\eta_{mf}=10^{3}$, $\phi_f = 0.2$ and consider test cases with varying contrast, $k_f = 10^3, 10^6$ and $10^6$ with $\kappa_f = k_f \cdot 10^{-20}$ [m$^2$]  \cite{akkutlu2017multiscale, akkutlu2018multiscale, talonov2016numerical, akkutlu2015multiscale}.

Parameters are set as follows \cite{fathi2014multi, akkutlu2018multiscale, akkutlu2015multiscale}
\begin{itemize}
\item $R = 8.31$[J/(K $\cdot$ mol)], $T = 323$[K], $Z = 1$,
\item $p_i = 20 \cdot  10^6$[Pa], $p_{w} = 5 \cdot  10^6$[Pa], $p_L = 10^6$[Pa],
\item $c_{init} = p_i/(ZRT)$[mol/m$^3$], $c_{w} = p_{well}/(ZRT)$[mol/m$^3$],
\item $\phi = 0.02$, $\phi_f = 0.2$, $\varepsilon_{ks} = \varepsilon_{kp} = 0.5$,
\item $D_s = D_i = D_k = 10^{-8}$[m$^2$/s],
\item $c_{\mu s} = 0.25 \cdot 10^5$[m$^2$/s], $K = Z R T/p_L$[m$^3$/mol],
\item $\mu = 10^{-5}$[Pa $\cdot$ s].
\end{itemize}
As initial condition, we set $c_{m,0} = c_{f,0} = c_{init}$ and  perform simulations for $T_{max} = 1$ day for Test 1 and  $T_{max} = 5$ days for Test 2 with $N_t=11$ time steps, $\tau = T_{max}/(N_T-1)$.
For a source term, we set $f_m=0$ and $f_f(c_f) = c_w ZRT \kappa_w (c_w - c_f)/\mu$ with $\kappa_w = k_w \cdot 10^{-20}$[m$^2$] and $k_w = 10^5$.

On the fine grid, we approximate the gas transport using first-order polynomial functions for both porous matrix and fractures. 
We consider two uniform coarse grids
\begin{itemize}
    \item $10 \times 10$ coarse grid with 121 nodes and 100 square cells;
    \item $20 \times 20$ coarse grid with 441 nodes and 400 square cells. 
\end{itemize}
The system sizes are $N_h = 64771$ for Test 1 ($N^m_h = 62975$ for porous matrix and $N^f_h = 1796$ for fractures)  and $DOF_h = 68110$ for Test 2 ($N^m_h = 64629$ for porous matrix and $N^f_h = 3481$ for fractures).  Implementation is performed using a Python programming language with Scipy sparse library to represent and solve the resulting system of linear equations at each time step. We use a PCG (preconditioned conjugate gradient) iterative solver with a proposed two-grid adaptive preconditioner for shale gas transport in fractured porous media. As a smoother, we use five Gauss-Seidel iterations ($\nu=5$) from pyamg library \cite{bell2022pyamg, bell2023pyamg}. Simulations are performed on a MacBook Pro with an Apple M2 Max chip (32 GB memory). 

\begin{figure}[h!]
\centering
\begin{subfigure}{0.24\textwidth}
\includegraphics[width=1 \textwidth]{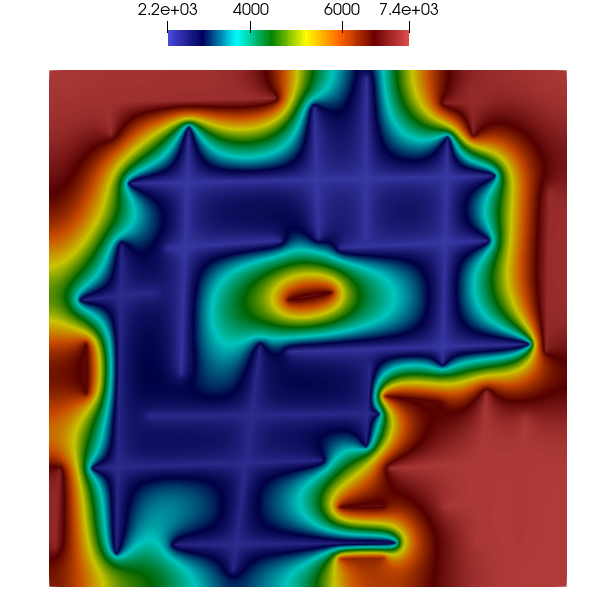}
\caption{Test 1a.}
\end{subfigure}
\begin{subfigure}{0.24\textwidth}
\includegraphics[width=1 \textwidth]{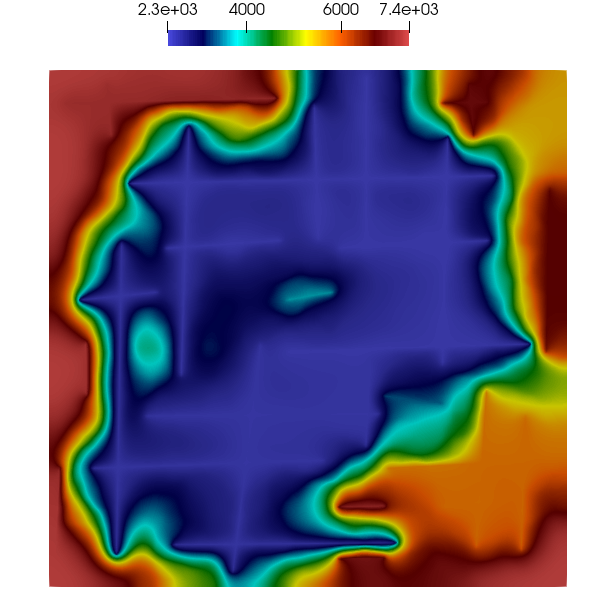}
\caption{Test 1b.}
\end{subfigure}
\begin{subfigure}{0.24\textwidth}
\includegraphics[width=1 \textwidth]{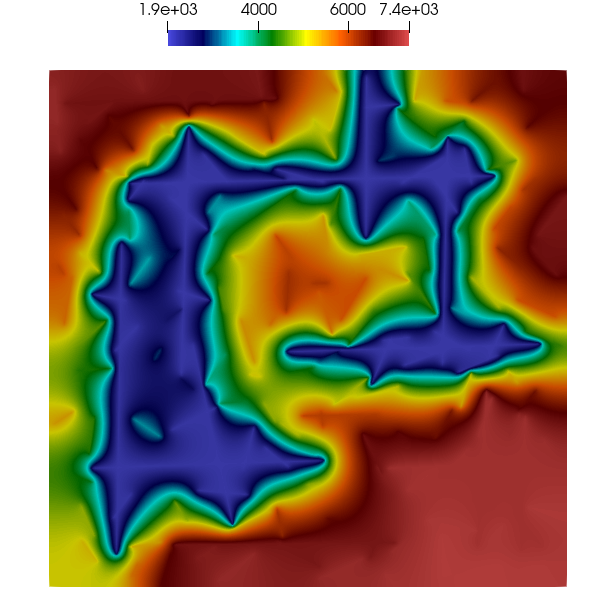}
\caption{Test 2a.}
\end{subfigure}
\begin{subfigure}{0.24\textwidth}
\includegraphics[width=1 \textwidth]{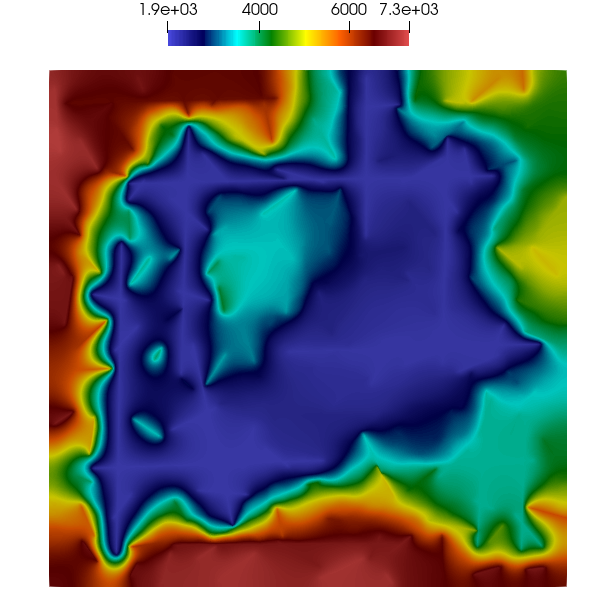}
\caption{Test 2b.}
\end{subfigure}
\caption{Solution at final time for homogeneous (first and third plots) and heterogeneous porous matrix (second and fourth plots) with contrast $10^9$.}
\label{fig:sol}
\end{figure}

In Figure \ref{fig:sol}, we plot a solution at the final time for $k_f = 10^9$. We  solve the system using an adaptive spectral preconditioner. In the first and second plots, we plot results for Test 1a and 1b. The results for Test 2a and Test 2b are presented in the third and fourth plots.

\begin{figure}[h!]
\centering
\includegraphics[width=1 \textwidth]{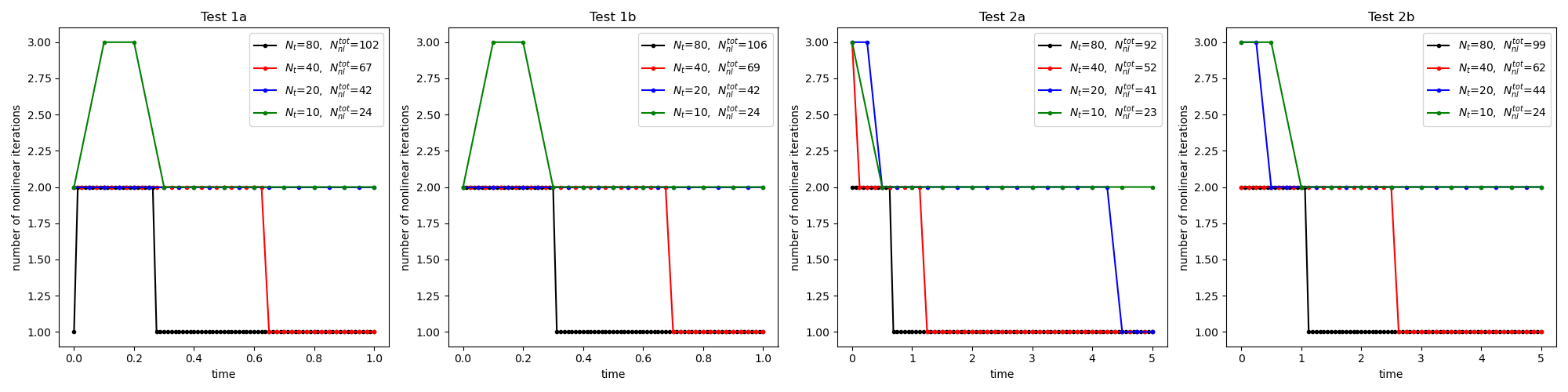}
\caption{Number of nonlinear iterations for implicit scheme with Picard iterations. Contrast $10^9$}
\label{fig:nonlin1}
\end{figure}

\begin{table}[h!]
\centering
\begin{tabular}{|c|cccc|}
\hline 
$N_t$ & Test 1a & Test 1b & Test 2a & Test 2b \\
\hline
\multicolumn{5}{|c|}{Contrast $10^3$}\\
\hline 
10 & 22 & 22 & 24 & 24 \\ 
20 & 42 & 42 & 41 & 44 \\ 
40 & 66 & 68 & 54 & 65 \\ 
80 & 94 & 99 & 93 & 100 \\ 
\hline
\multicolumn{5}{|c|}{Contrast $10^6$}\\
\hline 
10 & 24 & 24 & 23 & 24 \\ 
20 & 42 & 42 & 41 & 44 \\ 
40 & 67 & 69 & 52 & 62 \\ 
80 & 102 & 106 & 92 & 99 \\ 
\hline
\multicolumn{5}{|c|}{Contrast $10^9$}\\
\hline 
10 & 24 & 24 & 23 & 24 \\ 
20 & 42 & 42 & 41 & 44 \\ 
40 & 67 & 69 & 52 & 62 \\ 
80 & 102 & 106 & 92 & 99 \\ 
\hline 
\end{tabular}
\caption{Total number of nonlinear iterations for implicit scheme with Picard iterations}
\label{table-im-errf}
\end{table}

\begin{figure}[h!]
\centering
\includegraphics[width=1 \textwidth]{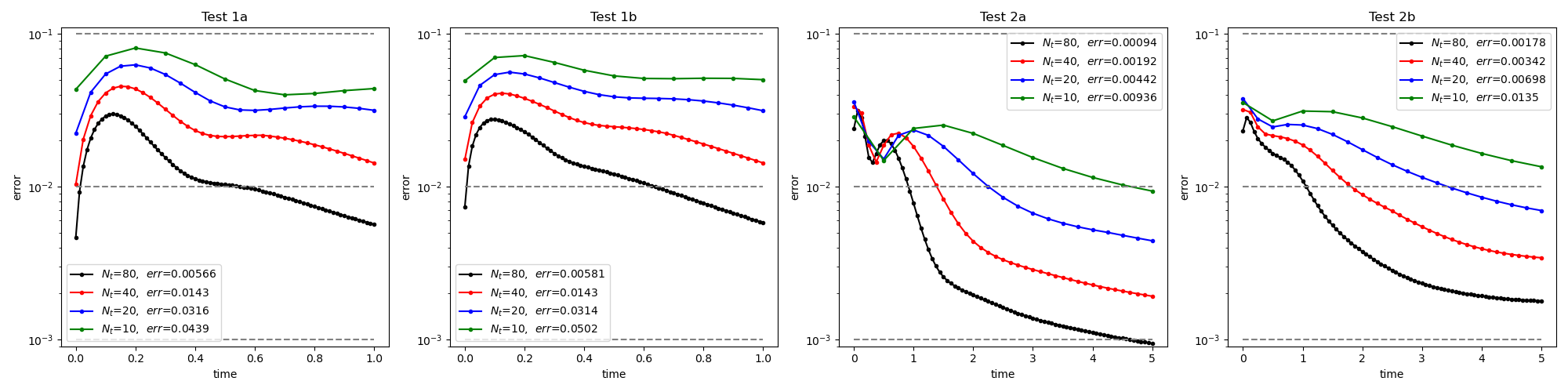}
\caption{Relative error dynamics for linearly implicit scheme. Contrast $10^9$}
\label{fig:nonlin2}
\end{figure}

\begin{table}[h!]
\centering
\begin{tabular}{|c|cccc|}
\hline 
$N_t$ & Test 1a & Test 1b & Test 2a & Test 2b \\
\hline
\multicolumn{5}{|c|}{Contrast $10^3$}\\
\hline 
10 & 0.03986 & 0.04450 & 0.00928 & 0.01367 \\ 
20 & 0.02914 & 0.02800 & 0.00427 & 0.00708 \\ 
40 & 0.01268 & 0.01244 & 0.00184 & 0.00346 \\ 
80 & 0.00462 & 0.00490 & 0.00086 & 0.00189 \\ 
\hline
\multicolumn{5}{|c|}{Contrast $10^6$}\\
\hline 
10 & 0.04392 & 0.05021 & 0.00936 & 0.01350 \\ 
20 & 0.03164 & 0.03143 & 0.00443 & 0.00699 \\ 
40 & 0.01430 & 0.01430 & 0.00192 & 0.00342 \\ 
80 & 0.00566 & 0.00581 & 0.00095 & 0.00178 \\ 
\hline
\multicolumn{5}{|c|}{Contrast $10^9$}\\
\hline 
10 & 0.04392 & 0.05022 & 0.00936 & 0.01350 \\ 
20 & 0.03164 & 0.03144 & 0.00442 & 0.00698 \\ 
40 & 0.01430 & 0.01430 & 0.00192 & 0.00342 \\ 
80 & 0.00566 & 0.00581 & 0.00094 & 0.00178 \\ 
\hline 
\end{tabular}
\caption{Relative error dynamics for linearly implicit scheme}
\label{table-ein-errf}
\end{table}

\begin{figure}[h!]
\centering
\includegraphics[width=0.8\linewidth]{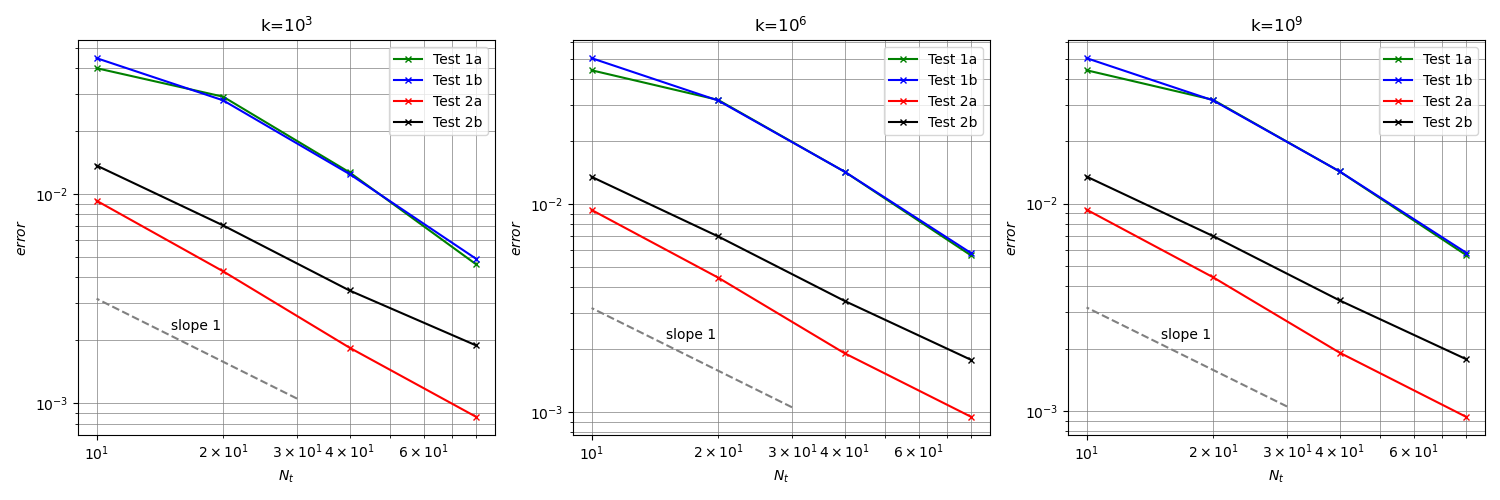}
\caption{Number of time steps ($N_t$) vs relative error at final time. Contrast $10^9$}
\label{fig:nonlin3}
\end{figure}

To compare two approaches for linearization, we calculate the relative error in $L_2$ norm
\[
e (c_{ref}, c) = \frac{||c_{ref} - c||_{L_2}}{||c_{ref}||_{L_2}} \times100 \%, \quad 
||c||^2_{L_2} = (c, c),
\]
where $c_{ref}$ is the reference solution using the implicit scheme and Picard's iterations for nonlinearity and $c$ is the solution using linearly implicit scheme. 
In nonlinear Picard's iterations, we set $N_{nl} = 10$ as a maximum number of nonlinear iterations and iterate till $e(c^{n+1, (m+1)}_h, c^{n+1, (m)}_h)\leq 0.1$ \%.
In Figure \ref{fig:nonlin1}, we plot the dynamics of the number of nonlinear iterations for an implicit scheme with Picard's iterations for contrast $10^9$. In Table \ref{table-im-errf}, we show a total number of nonlinear iterations $N_{nl}^{tot}$ for all three contrasts. 
We investigate the influence of the number of time steps $N_t = 10, 20, 40$, and $80$ on the nonlinear iterations. For a larger $N_t$, we obtain a faster decline of a number of nonlinear iterations. In the plot, we also show a total number of nonlinear iterations ($N^{tot}_{nl}$).

Next, we compare the proposed linearly implicit scheme with an implicit scheme.  In Table \ref{table-ein-errf}, we present relative $L_2$ error in percentage at a final time $e(c_{ref}, c)$,  where $c_{ref}$ denotes the reference solution using the implicit scheme and Picard's iterations for nonlinearity.
In Figure \ref{fig:nonlin2}, we plot the dynamics of the errors for the linearly implicit scheme for contrast $10^9$. 
In Figure \ref{fig:nonlin3}, we represent errors for contrast $10^9$ at a final time for different $N_t$ , showing a first-order convergence of the scheme.

\begin{table}[h!]
\centering
\begin{minipage}[t]{0.49\textwidth}
\centering
\begin{tabular}{| c | c | c | c |}
\hline
Contrast & $10^3$ & $10^6$ & $10^9$ \\
\hline 
$m$ ($N_H$) & \multicolumn{3}{|c|}{Test 1a} \\
\hline
1 (121) & $>$100 & $>$100 & $>$100 \\ 
2 (242) & 19.5 & 46.5 & $>$100 \\ 
4 (484) & 9.0 & 10.0 & 13.0 \\ 
8 (968) & 7.0 & 8.0 & 10.0 \\ 
12 (1452) & 6.0 & 7.0 & 8.0 \\ 
16 (1936) & 5.0 & 6.0 & 7.0 \\ 
\hline
$\delta_{\lambda} (N_H)$& \multicolumn{3}{|c|}{$m_{i}$ for $\lambda_{m_i} < \delta_{\lambda}$} \\
\hline
$10^{-4}$ (226) & 15.9 & 17.4 & 19.9 \\ 
$10^{-3}$ (285) & 14.0 & 14.7 & 17.0 \\ 
$10^{-2}$ (1293) & 6.0 & 7.0 & 8.6 \\ 
\hline
$\delta_{\lambda} (N_H)$& \multicolumn{3}{|c|}{$m_{i}+1$ for $\lambda_{m_i} < \delta_{\lambda}$} \\
\hline
$10^{-4}$ (347) & 12.0 & 13.0 & 15.0 \\ 
$10^{-3}$ (406) & 11.0 & 13.0 & 14.1 \\ 
$10^{-2}$ (1413) & 5.0 & 6.0 & 8.0 \\ 
\hline
\hline
$m$ ($N_H$)& \multicolumn{3}{|c|}{Test 1b } \\
\hline
1 (121) & $>$100 & $>$100 & $>$100 \\ 
2 (242) & 25.8 & $>$100 & $>$100 \\ 
4 (484) & 12.0 & 15.7 & 21.3 \\ 
8 (968) & 8.0 & 9.9 & 12.7 \\ 
12 (1452) & 6.0 & 8.0 & 9.9 \\ 
16 (1936) & 5.0 & 7.0 & 8.0 \\ 
\hline 
$\delta_{\lambda} (N_H)$& \multicolumn{3}{|c|}{$m_{i}$ for $\lambda_{m_i} < \delta_{\lambda}$} \\
\hline
$10^{-4}$ (228) & 19.9 & 28.2 & 34.8 \\ 
$10^{-3}$ (291) & 18.0 & 23.2 & 28.4 \\ 
$10^{-2}$ (1288) & 6.0 & 8.0 & 9.0 \\ 
\hline
$\delta_{\lambda} (N_H)$& \multicolumn{3}{|c|}{$m_{i}+1$ for $\lambda_{m_i} < \delta_{\lambda}$} \\
\hline
$10^{-4}$ (349) & 13.9 & 17.7 & 22.5 \\ 
$10^{-3}$ (412) & 13.0 & 16.0 & 21.1 \\ 
$10^{-2}$ (1409) & 6.0 & 8.0 & 9.0 \\ 
\hline
\end{tabular}
\end{minipage}
\ \ 
\begin{minipage}[t]{0.49\textwidth}
\centering
\begin{tabular}{| c | c | c | c |}
\hline
Contrast & $10^3$ & $10^6$ & $10^9$ \\
\hline 
$m$ ($N_H$) & \multicolumn{3}{|c|}{Test 2a} \\
\hline
1 (121) & $>$100 & $>$100 & $>$100 \\ 
2 (242) & $>$100 & $>$100 & $>$100 \\ 
4 (484) & $>$100 & $>$100 & $>$100 \\ 
8 (968) & 22.4 & $>$100 & $>$100 \\ 
12 (1452) & 8.3 & 10.4 & 14.1 \\ 
16 (1936) & 7.1 & 9.0 & 12.0 \\ 
\hline
$\delta_{\lambda} (N_H)$& \multicolumn{3}{|c|}{$m_{i}$ for $\lambda_{m_i} < \delta_{\lambda}$} \\
\hline
$10^{-4}$ (802) & 18.9 & 20.4 & 24.9 \\ 
$10^{-3}$ (821) & 14.9 & 17.2 & 22.0 \\ 
$10^{-2}$ (1631) & 7.4 & 9.0 & 12.0 \\ 
\hline
$\delta_{\lambda} (N_H)$& \multicolumn{3}{|c|}{$m_{i}+1$ for $\lambda_{m_i} < \delta_{\lambda}$} \\
\hline
$10^{-4}$ (923) & 13.0 & 15.5 & 20.1 \\ 
$10^{-3}$ (942) & 12.7 & 15.5 & 20.0 \\ 
$10^{-2}$ (1752) & 6.3 & 8.0 & 11.0 \\ 
\hline
\hline
$m$ ($N_H$)& \multicolumn{3}{|c|}{Test 2b} \\
\hline
1 (121) & $>$100 & $>$100 & $>$100 \\ 
2 (242) & $>$100 & $>$100 & $>$100 \\ 
4 (484) & $>$100 & $>$100 & $>$100 \\ 
8 (968) & 20.0 & $>$100 & $>$100 \\ 
12 (1452) & 8.9 & 11.5 & 15.2 \\ 
16 (1936) & 7.2 & 9.0 & 11.7 \\ 
\hline
$\delta_{\lambda} (N_H)$& \multicolumn{3}{|c|}{$m_{i}$ for $\lambda_{m_i} < \delta_{\lambda}$} \\
\hline
$10^{-4}$ (802) & 20.2 & 23.3 & 26.3 \\ 
$10^{-3}$ (821) & 15.4 & 17.6 & 23.5 \\ 
$10^{-2}$ (1644) & 7.3 & 9.0 & 11.5 \\ 
\hline
$\delta_{\lambda} (N_H)$& \multicolumn{3}{|c|}{$m_{i}+1$ for $\lambda_{m_i} < \delta_{\lambda}$} \\
\hline
$10^{-4}$ (923) & 14.4 & 16.8 & 21.5 \\ 
$10^{-3}$ (942) & 13.7 & 15.6 & 20.5 \\ 
$10^{-2}$ (1765) & 7.0 & 8.5 & 11.3 \\ 
\hline
\end{tabular}
\end{minipage}
\caption{Average number of iterations for different numbers of multiscale basis functions. $10 \times 10$ coarse grid and $\nu=5$ smoothing iterations. $N_H = m \cdot 121$ (m basis functions), $N_H = \sum_i m_{i}$ (adaptive coarse space), $N_h = 64771$ (Test 1a and 1b) and $N_h = 68110$ (Test 2a and 2b). }
\label{table-ms-10}
\end{table}

\begin{table}[h!]
\centering
\begin{minipage}[t]{0.49\textwidth}
\centering
\begin{tabular}{| c | c | c | c |}
\hline
Contrast & $10^3$ & $10^6$ & $10^9$ \\
\hline 
$m$ ($N_H$) & \multicolumn{3}{|c|}{Test 1a} \\
\hline
1 (441) & 17.6 & 40.1 & $>$100 \\ 
2 (882) & 9.0 & 10.0 & 12.7 \\ 
4 (1764) & 5.5 & 6.0 & 7.0 \\ 
8 (3528) & 4.0 & 5.0 & 6.0 \\ 
12 (5292) & 3.0 & 4.0 & 5.0 \\ 
16 (7056) & 3.0 & 4.0 & 4.0 \\ 
\hline
$\delta_{\lambda} (N_H)$& \multicolumn{3}{|c|}{$m_{i}$ for $\lambda_{m_i} < \delta_{\lambda}$} \\
\hline
$10^{-4}$ (486) & 10.0 & 11.7 & 14.2 \\ 
$10^{-3}$ (531) & 10.0 & 11.4 & 14.0 \\ 
$10^{-2}$ (1629) & 6.0 & 6.0 & 8.0 \\ 
\hline
$\delta_{\lambda} (N_H)$& \multicolumn{3}{|c|}{$m_{i}+1$ for $\lambda_{m_i} < \delta_{\lambda}$} \\
\hline
$10^{-4}$ (927) & 9.0 & 9.0 & 10.0 \\ 
$10^{-3}$ (972) & 9.0 & 9.0 & 10.0 \\ 
$10^{-2}$ (2070) & 5.0 & 6.0 & 7.0 \\ 
\hline
\hline
$m$ ($N_H$)& \multicolumn{3}{|c|}{Test 1b } \\
\hline
1 (441) & 22.5 & $>$100 & $>$100 \\ 
2 (882) & 9.7 & 14.9 & 18.9 \\ 
4 (1764) & 6.0 & 7.0 & 9.5 \\ 
8 (3528) & 4.0 & 5.0 & 6.0 \\ 
12 (5292) & 3.0 & 4.0 & 5.0 \\ 
16 (7056) & 3.0 & 4.0 & 4.6 \\ 
\hline 
$\delta_{\lambda} (N_H)$& \multicolumn{3}{|c|}{$m_{i}$ for $\lambda_{m_i} < \delta_{\lambda}$} \\
\hline
$10^{-4}$ (486) & 16.5 & 21.9 & 28.7 \\ 
$10^{-3}$ (558) & 13.0 & 18.9 & 26.0 \\ 
$10^{-2}$ (1623) & 6.5 & 8.0 & 10.0 \\ 
\hline
$\delta_{\lambda} (N_H)$& \multicolumn{3}{|c|}{$m_{i}+1$ for $\lambda_{m_i} < \delta_{\lambda}$} \\
\hline
$10^{-4}$ (927) & 9.0 & 12.0 & 14.5 \\ 
$10^{-3}$ (999) & 9.0 & 11.6 & 14.0 \\ 
$10^{-2}$ (2064) & 5.9 & 7.0 & 9.8 \\ 
\hline
\end{tabular}
\end{minipage}
\ \ 
\begin{minipage}[t]{0.49\textwidth}
\centering
\begin{tabular}{| c | c | c | c |}
\hline
Contrast & $10^3$ & $10^6$ & $10^9$ \\
\hline 
$m$ ($N_H$) & \multicolumn{3}{|c|}{Test 2a} \\
\hline
1 (441) & $>$100 & $>$100 & $>$100 \\ 
2 (882) & $>$100 & $>$100 & $>$100 \\ 
4 (1764) & 12.5 & 28.0 & 72.0 \\ 
8 (3528) & 5.0 & 6.0 & 8.0 \\ 
12 (5292) & 4.0 & 5.0 & 6.0 \\ 
16 (7056) & 3.0 & 4.0 & 5.0 \\ 
\hline
$\delta_{\lambda} (N_H)$& \multicolumn{3}{|c|}{$m_{i}$ for $\lambda_{m_i} < \delta_{\lambda}$} \\
\hline
$10^{-4}$ (1187) & 14.5 & 16.7 & 21.0 \\ 
$10^{-3}$ (1204) & 14.3 & 16.4 & 20.6 \\ 
$10^{-2}$ (2163) & 7.0 & 8.0 & 10.9 \\ 
\hline
$\delta_{\lambda} (N_H)$& \multicolumn{3}{|c|}{$m_{i}+1$ for $\lambda_{m_i} < \delta_{\lambda}$} \\
\hline
$10^{-4}$ (1628) & 9.4 & 10.4 & 13.0 \\ 
$10^{-3}$ (1645) & 9.4 & 10.4 & 13.0 \\ 
$10^{-2}$ (2605) & 6.0 & 7.0 & 9.9 \\ 
\hline
\hline
$m$ ($N_H$)& \multicolumn{3}{|c|}{Test 2b} \\
\hline
1 (441) & $>$100 & $>$100 & $>$100 \\ 
2 (882) & $>$100 & $>$100 & $>$100 \\ 
4 (1764) & 13.5 & 30.5 & 86.0 \\ 
8 (3528) & 5.2 & 7.0 & 8.0 \\ 
12 (5292) & 4.0 & 4.5 & 6.0 \\ 
16 (7056) & 3.0 & 4.0 & 5.0 \\ 
\hline
$\delta_{\lambda} (N_H)$& \multicolumn{3}{|c|}{$m_{i}$ for $\lambda_{m_i} < \delta_{\lambda}$} \\
\hline
$10^{-4}$ (1187) & 19.0 & 18.5 & 24.5 \\ 
$10^{-3}$ (1212) & 15.6 & 18.4 & 23.0 \\ 
$10^{-2}$ (2149) & 7.3 & 8.4 & 11.0 \\ 
\hline
$\delta_{\lambda} (N_H)$& \multicolumn{3}{|c|}{$m_{i}+1$ for $\lambda_{m_i} < \delta_{\lambda}$} \\
\hline
$10^{-4}$ (1628) & 9.7 & 11.0 & 13.7 \\ 
$10^{-3}$ (1653) & 9.6 & 11.0 & 13.7 \\ 
$10^{-2}$ (2590) & 6.0 & 7.5 & 10.0 \\ 
\hline
\end{tabular}
\end{minipage}
\caption{Average number of iterations for different numbers of multiscale basis functions. $20 \times 20$ coarse grid and $\nu=5$ smoothing iterations. $N_H = m \cdot 441$ (m basis functions), $N_H = \sum_i m_{i}$ (adaptive coarse space), $N_h = 64771$ (Test 1a and 1b) and $N_h = 68110$ (Test 2a and 2b). }
\label{table-ms-20}
\end{table}

\begin{figure}[h!]
\centering
\begin{subfigure}{1\textwidth}
\includegraphics[width=0.23 \textwidth]{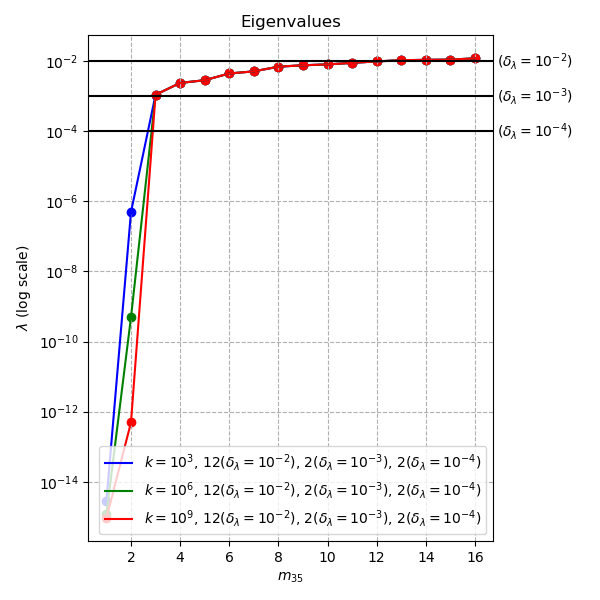}
\includegraphics[width=0.24 \textwidth]{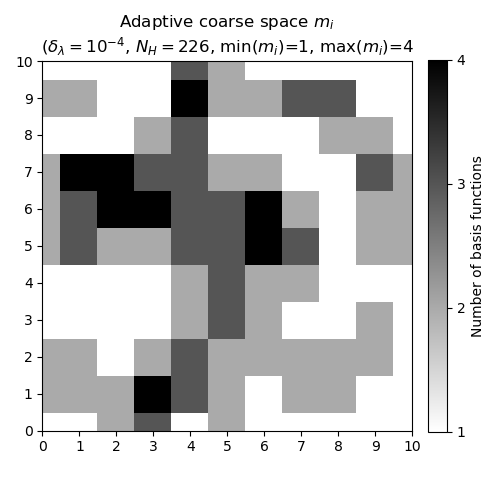}
\includegraphics[width=0.24 \textwidth]{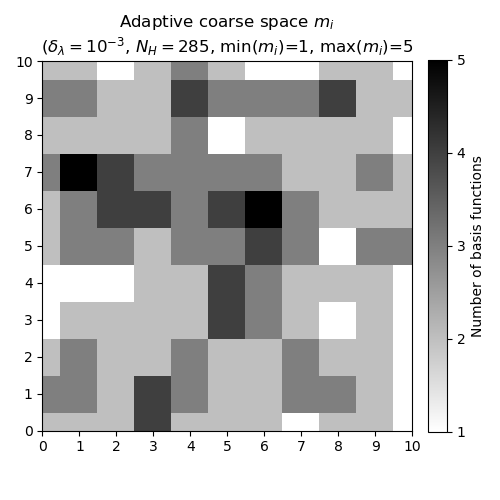}
\includegraphics[width=0.24 \textwidth]{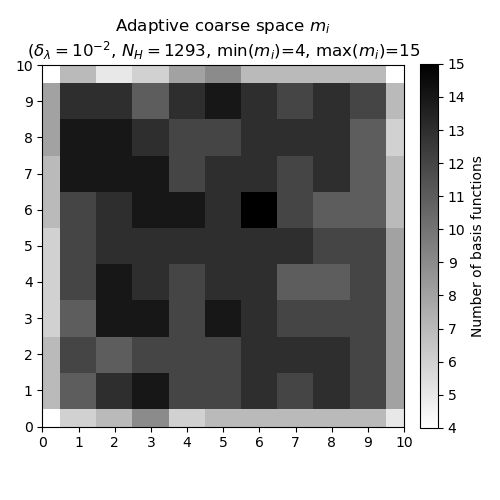}
\caption{Test 1a. $10 \times 10$ coarse grid.}
\end{subfigure}
\begin{subfigure}{1\textwidth}
\includegraphics[width=0.24 \textwidth]{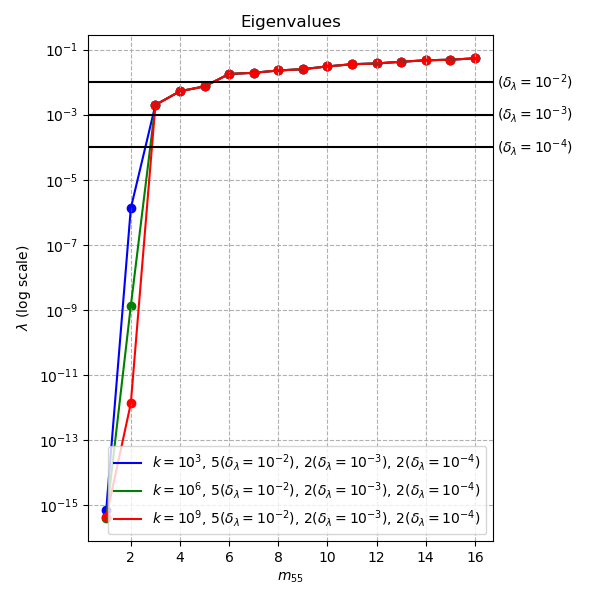}
\includegraphics[width=0.24 \textwidth]{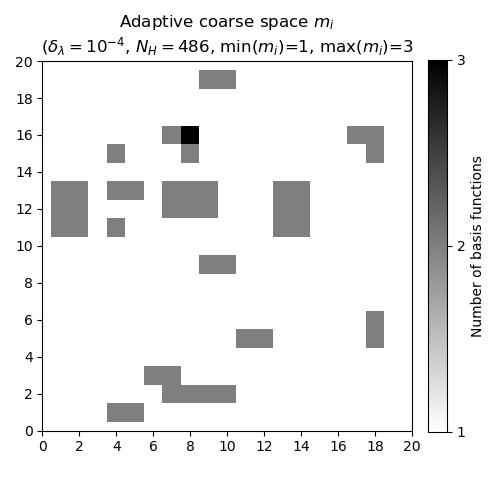}
\includegraphics[width=0.24 \textwidth]{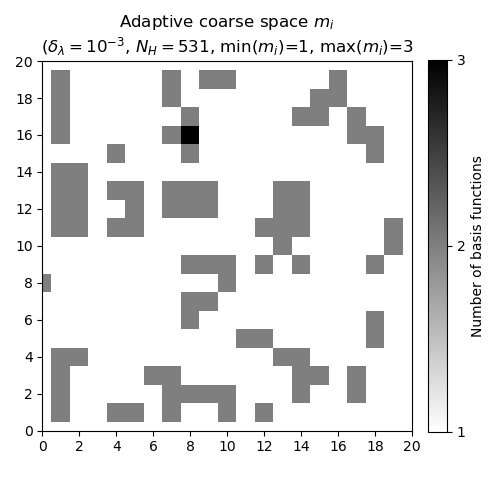}
\includegraphics[width=0.24 \textwidth]{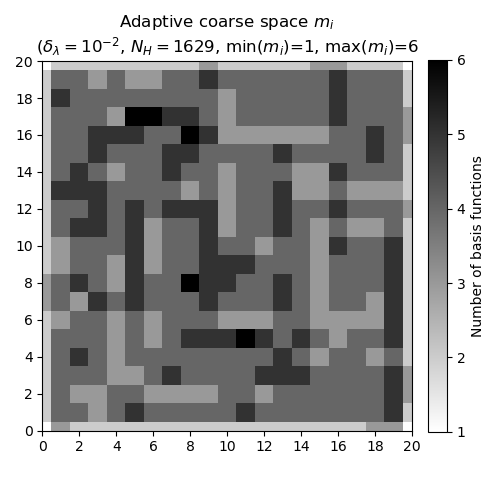}
\caption{Test 1a. $20 \times 20$ coarse grid.}
\end{subfigure}
\begin{subfigure}{1\textwidth}
\includegraphics[width=0.24 \textwidth]{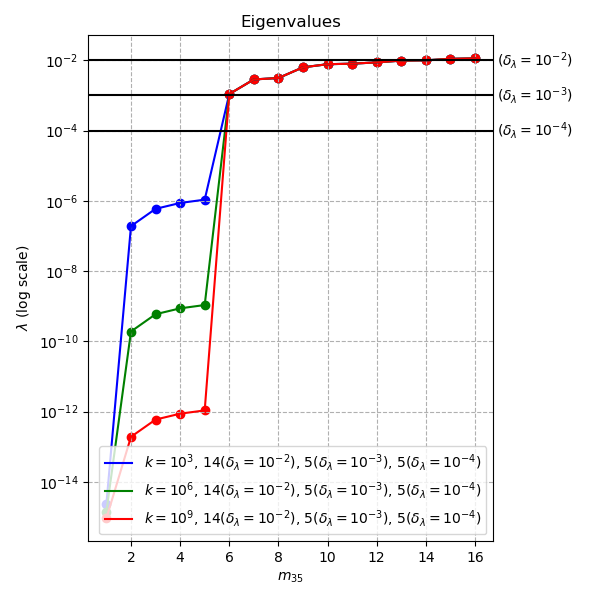}
\includegraphics[width=0.24 \textwidth]{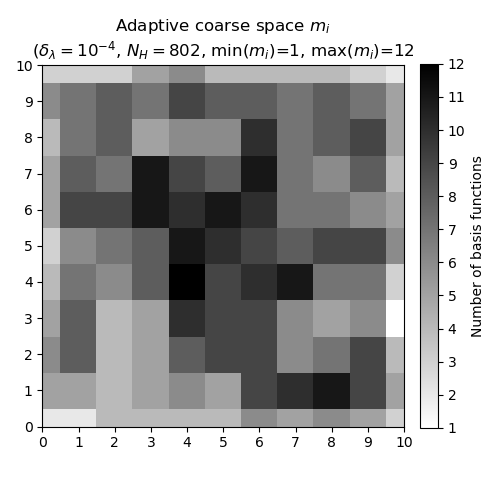}
\includegraphics[width=0.24 \textwidth]{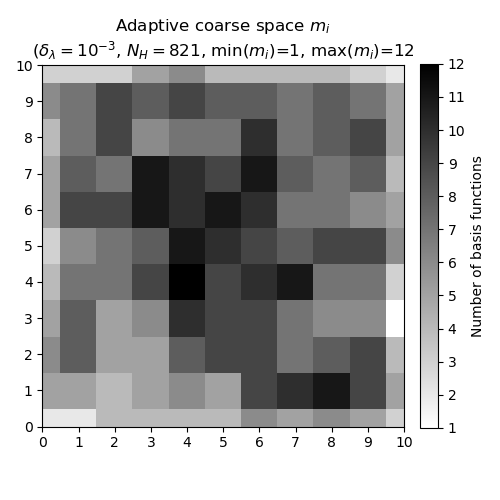}
\includegraphics[width=0.24 \textwidth]{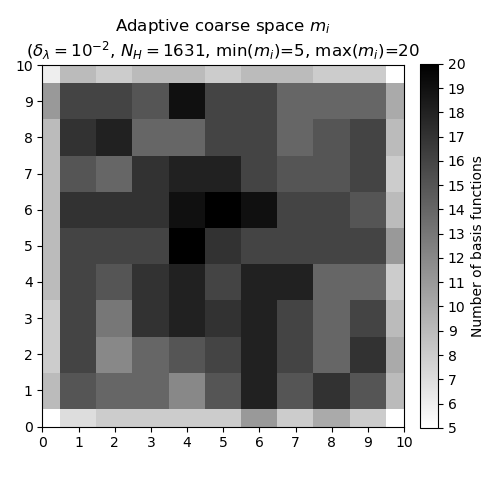}
\caption{Test 2a. $10 \times 10$ coarse grid.}
\end{subfigure}
\begin{subfigure}{1\textwidth}
\includegraphics[width=0.24 \textwidth]{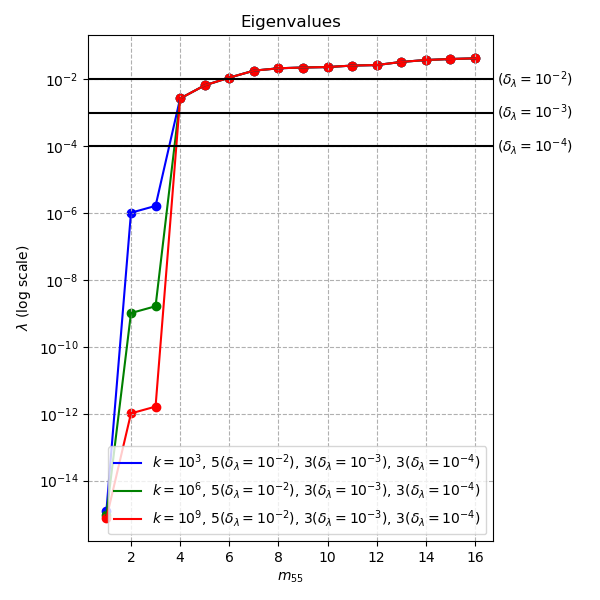}
\includegraphics[width=0.24 \textwidth]{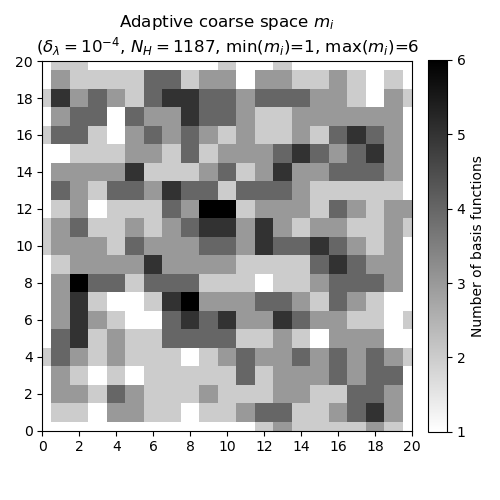}
\includegraphics[width=0.24 \textwidth]{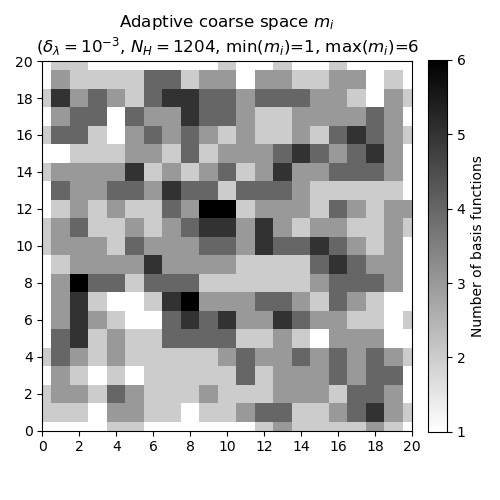}
\includegraphics[width=0.24 \textwidth]{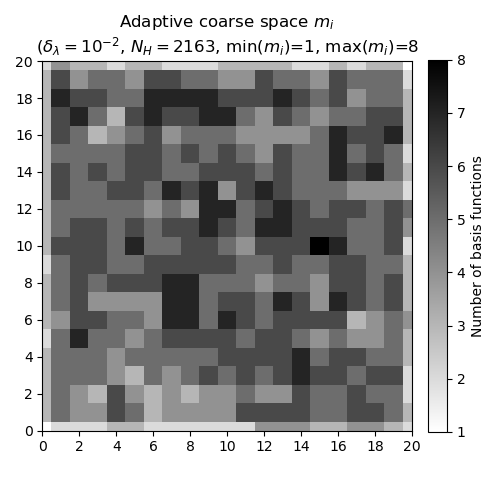}
\caption{Test 2a. $20 \times 20$ coarse grid.}
\end{subfigure}
\caption{Illustration of adaptivity in spectral space construction for Test 1a and 2a.   First column: first smallest eigenvalues $\omega_{35}$ on $10 \times 10$ coarse grid and $\omega_{55}$ on $20 \times 20$ coarse grid. 
Second, third and fourth columns: adaptivity for $\delta_{\lambda} = 10^{-4}, 10^{-3}$ and $10^{-2}$ (from left to right) with $k_f=10^9$.}
\label{fig:adp}
\end{figure}

In Tables \ref{table-ms-10} and \ref{table-ms-20}, we present results for an iterative method with a multiscale preconditioner using $10 \times 10$ and $20 \times 20$ coarse grids. The influence of permeability contrast is investigated for $10^3, 10^6$, and $10^9$ contrasts. 
The study uses 100 in an iterative solver as the maximum number of iterations and relative residual tolerance $\varepsilon =10^{-9}$ for iterative solver. The convergence results are presented for different fixed numbers of multiscale basis functions $m\in\{1,2,4,6,8,12,16\}$ as well as for an adaptive choice of multiscale basis functions, where we choose a different number of multiscale basis functions in each local domain by the given threshold of eigenvalues, $m_i =\max i$ such that $\lambda_{m_i} < \delta_{\lambda}$ with $\delta_{\lambda} \in\{10^{-2}, 10^{-3}, 10^{-4}\}$. 
Moreover, we consider two cases of adaptivity, $m_{i}$ local basis functions and $m_i + 1$ functions, leading to  $N_H = \sum_{i=1}^N m_{i}$ and  $N_H = \sum_{i=1}^N (m_{i}+1) = \sum_{i=1}^N m_{i} + N$, respectively. 
In the first column, we specify the number or adaptive, tolerance $\delta_{\lambda}$ for coarse basis functions and the size of the resulting coarse grid system. In the proceeding columns to the right, we show the average number of iterations in each time step ($N^{it}_{Av}$).

First, we emphasize that for all four test problems, two coarse-grid sizes and a wide range of permeability contrasts, with $\mathcal{O}(1)$ local multiscale basis functions per subdomain we are able to achieve excellent two-grid convergence. Indeed, with sufficient local basis functions, iteration counts are robust across all problem formulations tested. Given a residual tolerance of $10^{-9}$, observed low iteration counts correspond to average convergence factors of $0.1$ or less. The key is having sufficient local basis functions, as without a sufficiently accurate coarse-grid operator convergence deteriorates quickly. For example, for a $10 \times 10$ coarse grid, the iterative method does not converge in 100 iterations for $m=2$ in Test 1a, 2a, and 2b for high contrast ($10^9$). For smaller contrast, we can use two basis functions in Test 1a and 1b. In Test 1b, with 30 fractures and a heterogeneous background, we obtain convergence with two bases; however, using four basis functions can obtain a two times smaller number of iterations, and the number of iterations is further reduced with better coarse scale approximation. Only 5-8 iterations are required for $m=16$ for Tests 1a and 1b. Similarly, for the geometry with a larger number of fractures in Test 2a and 2b, we observe that the iterative solver does not converge in 100 iterations when we use $m=1,2$ and $4$ basis functions per local domain. However, we see that when we take a sufficient number of basis functions, we obtain a solution in 7-12 iterations for both Test 2a and 2b. Tables \ref{table-ms-20} show the number of iterations on $20 \times 20$ coarse grid and have a similar dependence on a number of local basis functions. 

To that end, the adaptive approach proves to be both robust and efficient in all cases, achieving excellent convergence for a minimal number of basis functions. For example, Test 2 has many more fractures than Test 1, and typically requires $2-3\times$ as many fixed local basis functions $m$ to achieve fast convergence in $\lessapprox 12$ iterations (across different permeability contrasts and coarse-grid sizes). The adaptive approach is quite robust in this respect however, for a fixed tolerance, e.g. $\delta_{\lambda} = 10^{-3}$, achieving fairly similar iteration counts across \emph{all} problem formulations. We note that adding one additional basis function per local domain in the adaptive approach does not notably improve the number of iterations but raises the number of degrees of freedom on the coarse grid. As in most adaptive schemes, one must make a choice in parameter $\delta_{\lambda}$, but here we demonstrate that one choice remains robust across many different problems. In general, $\delta_{\lambda} = 10^{-2}$ seems to be overkill in most cases, increasing coarse-grid size and number of local basis functions more than the benefit in improvement in iteration count, but $\delta_{\lambda} = 10^{-3}$ and $\delta_{\lambda} = 10^{-4}$ converge effectively in all cases. Example eigenvalue behavior for local subdomains of various problem formulations are shown in first column of Figure \ref{fig:adp}.

In addition, by setting a local tolerance rather than a fixed number of basis functions, we adaptively choose the number of basis functions depending on the local subdomain. The ``right'' choice can vary significantly across domain, as seen in Figures \ref{fig:adp} for Test 1a and 2a. Here, we plot how many basis functions are chosen across the domain for various problems, seeing as much variation in a single global problem from 1 local basis function to 12 local basis functions selected in different subdomains. The result is the robust convergence discussed previously at a reduced cost compared with a globally fixed number of basis functions. For example, in Test 2a with $20\times 20$ coarse grid (see Table \ref{table-ms-20}), adaptive iteration counts with $\delta_{\lambda}=10^{-3}$ match iteration counts for fixed $m=4$, while requiring 1204 coarse DOFs compared with 1764. For Test 2a with homogeneous background and $20\times 20$ coarse grid (see Table \ref{table-ms-20}), for $k_f=10^9$, $\delta_{\lambda}=10^{-3}$ converges in $3\times$ fewer iterations than $m=4$, while also having 50\% less coarse DOFs. Not all cases offer such notable improvements, but in all cases, we see reduced coarse grid size for a given iteration count when comparing adaptive to a fixed number of basis functions. Overall, we can state that the adaptive approach is a robust tool for constructing accurate multiscale coarse approximations and converging in a few iterations.

\section{Conclusion}

We consider a time-dependent mixed-dimensional model for shale gas transport in fractured porous media. To construct a discrete system, we use an explicit-implicit scheme and finite element discretization in space. 
A careful choice of an additive partitioning of the nonlinear operators is proposed to separate the linear and nonlinear parts and ensure the stability of the time approximation scheme, with fixed linear approximation of stiff dynamics for all time. To construct an efficient two-grid preconditioner, we construct an adaptive spectral multiscale space for shale gas transport in fractured porous media. The proposed adaptive approach automatically chooses the number of local multiscale basis functions across different subdomains, depending on the local character of the matrix. We show that the method gives an optimal convergence rate for fractured porous media, achieving nine orders of magnitude residual reduction in $\mathcal{O}(1)$ iterations. The multiple degrees of freedom in this spectral approach can be associated with the multicontinuum approach, where each basis represents continua in each local domain. Numerical results are presented to illustrate the convergence of the proposed method in a two-dimensional domain with two given fracture distributions (30 fractures and 160 fractures) for settings with homogeneous and heterogeneous shale matrix properties. The results show that we can obtain a contrast-independent convergence within the adaptive choice of continua necessary to represent a complex shale gas transport in fractured porous media. The linearly implicit time integration scheme demonstrates good stability, expected convergence order, and prevents the expensive recomputation of a multiscale coarse space every time step. In addition, the adaptive approach is demonstrated to be robust with a single parameter choice across all problem formulations considered. 

\section*{Acknowledgments}
BSS was supported by the DOE Office of Advanced Scientific Computing Research Applied Mathematics program through Contract No. 89233218CNA000001. Los Alamos National Laboratory report number LA-UR-25-25052.

\section*{Appendix. Proof of stability}\label{secA1}

In this section, we proof Theorem \ref{t:t1}. 
Noting the identities
\[
\begin{split}
&2(c^{n+1} - c^n) = (c^{n+1} - c^{n-1}) + (c^{n+1} - 2 c^n + c^{n-1}), \\ 
&2(c^{n} - c^{n-1}) = (c^{n+1} - c^{n-1}) - (c^{n+1} - 2 c^n + c^{n-1}), \\
&2 c^{n+1} = (c^{n+1} - c^{n-1}) + (c^{n+1} - 2 c^n + c^{n-1}) + 2 c^n,
\end{split}
\]
then we rewrite \eqref{ein} as follows
\[
\begin{split}
S^{(lin)} &\frac{c^{n+1} - c^n}{\tau} 
+ S^{(nl)}(c^n) \frac{c^{n} - c^{n-1}}{\tau} 
+ D^{(lin)} c^{n+1} + D^{(nl)}(c^n) c^n \\
&= \frac{1}{2 \tau}(S^{(lin)} + S^{(nl)}(c^n) + \tau D^{(lin)}) (c^{n+1} - c^{n-1})\\
&+ \frac{1}{2 \tau}(S^{(lin)} - S^{(nl)}(c^n) + \tau D^{(lin)}) (c^{n+1} - 2 c^n + c^{n-1})\\
&+ (D^{(lin)} +  D^{(nl)}(c^n)) c^n
= F.
\end{split}
\]
Therefore, we rewrite scheme \eqref{ein} in a unified way known as a canonical form of three-step schemes \cite{samarskii2001theory, samarskii2001additive, vabishchevich2020explicit}
\begin{equation}
\label{cform}
\begin{split}
B^n \frac{c^{n+1} - c^{n-1}}{2\tau} 
+ R^n (c^{n+1} - 2c^n + c^{n-1}) 
+  D^n c^n = F.
\end{split}
\end{equation}
with 
\[
D^n = D^{(lin)} + D^{(nl)}(c^n), \quad 
B^n = S^{(lin)} + S^{(nl)}(c^n) + \tau D^{(lin)} > 0, 
\]\[
R^n = \frac{1}{2 \tau} (S^{(lin)} - S^{(nl)} (c^n) + \tau D^{(lin)}) > 0.
\]

We let
\[
c^n + c^{n-1} = 2 y^n, \quad 
c^n - c^{n-1} = w^n,
\]
then 
\[
\begin{split}
&w^{n+1} + w^{n} = c^{n+1} - c^{n-1}, \\
&w^{n+1} - w^{n} = c^{n+1} - 2 c^n + c^{n-1} = 2(y^{n+1} - y^{n}),\\
&2(y^{n+1} + y^{n}) = c^{n+1} + 2 c^{n} + c^{n-1}, \\
& c^n 
= \frac{1}{4}(c^{n+1} + 2 c^{n} + c^{n-1}) - \frac{1}{4} (c^{n+1} - 2 c^{n} + c^{n-1})
= \frac{1}{2}(y^{n+1} + y^{n}) - \frac{1}{4}(w^{n+1} - w^{n}).
\end{split}
\]

Next, we multiply \eqref{cform} by $w^{n+1} + w^{n} = 2(y^{n+1} - y^{n})$ and obtain 
\begin{equation}
\label{t1b}
\begin{split}
\frac{1}{2\tau} &(B^n (c^{n+1} - c^{n-1}),w^{n+1} + w^{n}) 
\\
&+ (R^n (c^{n+1} - 2c^n + c^{n-1}), w^{n+1} + w^{n}) 
+ (D^n c^n, w^{n+1} + w^{n}) \\
& = 
\frac{1}{2\tau} (B^n(w^{n+1} + w^{n}),w^{n+1} + w^{n}) \\
&+ \left(\left(R^n - \frac{1}{4} D^n \right) (w^{n+1} - w^{n}), w^{n+1} + w^{n} \right) \\
&+ (D^n (y^{n+1} + y^{n}), y^{n+1} - y^{n}) 
= (F, w^{n+1} + w^{n}).
\end{split}
\end{equation}

With $D^n = (D^n)^T \geq 0$, $R^n = (R^n)^T$ and 
\[
\begin{split}
R^n - \frac{1}{4} D^n 
&= \frac{1}{2 \tau} (S^{(lin)} - S^{(nl)}(c^n) + \tau D^{(lin)})
- \frac{1}{4}  D^n\\
&= \frac{1}{2 \tau} (S^{(lin)} - S^{(nl)}(c^n)) + \frac{1}{4}(D^{(lin)} - D^{(nl)}(c^n)).
\end{split}
\]
By assumption, we have $R^n - \tfrac{1}{4}D^n > 0$, and thus
\[
\begin{split}
&\left(\left(R^n - \frac{1}{4} D^n \right) (w^{n+1} - w^{n}), w^{n+1} + w^{n}\right)\\
& \quad \quad  = \left(\left(R^n - \frac{1}{4} D^n \right) w^{n+1}, w^{n+1} \right) 
 - \left(\left(R^n - \frac{1}{4} D^n \right) w^{n}, w^{n} \right)
\end{split}
\]
and 
$(D^n (y^{n+1} - y^{n}), y^{n+1} + y^{n}) 
= (D^n y^{n+1}, y^{n+1}) - (D^n y^{n},y^{n})$.

By the Cauchy?Schwarz and Young's inequalities, we obtain
\begin{equation}
\label{csh}
(F, w^{n+1} + w^{n}) 
\leq 
\frac{1}{4 \varepsilon}  \|w^{n+1} + w^{n}\|_{X^n}^2 
+ \varepsilon \| F \|_{(X^n)^{-1}}^2.
\end{equation}
Then from \eqref{t1b} with $\varepsilon=\tau/2$, we have 
\[
\begin{split}
\frac{1}{2\tau} &((B^n-X^n) (w^{n+1} + w^{n}),w^{n+1} + w^{n} ) \\
&+ \left(\left(R^n - \frac{1}{4} D^n \right) w^{n+1}, w^{n+1}  \right) 
+ (D^n y^{n+1}, y^{n+1}) \\
&\leq  
\left(\left(R^n - \frac{1}{4} D^n \right) w^{n}, w^{n}  \right) 
+ (D^n y^{n}, y^{n}) 
+ \frac{\tau}{2} \| F \|_{(X^n)^{-1}}^2.
\end{split}
\]
Therefore for $X^n = B^n$ and $X^n = S^n$, we have
\[
\mathcal{E}^{n+1} \leq  \mathcal{E}^{n} + \frac{\tau}{2} \| F \|_{(B^n)^{-1}}^2
\quad \text{ or } \quad 
\mathcal{E}^{n+1} \leq  \mathcal{E}^{n} + \frac{\tau}{2} \| F \|_{(S^n)^{-1}}^2.
\]
with 
\[
\begin{split}
\mathcal{E}^{n} 
= \| w^{n}\|_{(R^n - \frac{1}{4} D^n)}^2
+ \|y^{n}\|_{D^n}^2 
= \|c^{n} - c^{n-1}\|_{(R^n - \frac{1}{4} D^n)}^2
+ \|c^{n} + c^{n-1}\|_{D^n}^2
\end{split}
\]

Finally, we have an unconditionally stable scheme \eqref{ein} with  
\[
R^n - \frac{1}{4} D^n 
= \frac{1}{2 \tau} (S^{(lin)} - S^{(nl)}(c^n)) + \frac{1}{4}(D^{(lin)} - D^{(nl)}(c^n)) 
> 0.
\]
or
\[
S^{(lin)} - S^{(nl)}(c^n) > 0 \quad 
D^{(lin)} - D^{(nl)}(c^n) > 0,
\]
for all $n=1,\ldots,N_t$.

\bibliographystyle{plain}
\bibliography{lit}

\end{document}